%% file: main.tex
\documentclass{article}

\usepackage[T1]{fontenc}
\usepackage[utf8]{inputenc}

\usepackage{amsmath}
\usepackage{amssymb}

\usepackage[
    backend=biber,
    style=numeric,
    sorting=none
]{biblatex}

\usepackage{graphicx}
\usepackage{subcaption}
\usepackage[percent]{overpic}
\usepackage{xcolor}

\usepackage{subcaption}

\usepackage{overpic}
\usepackage{subcaption}
\usepackage{relsize}
\usepackage{tikz}
\usetikzlibrary{calc}

\usepackage{hyperref}
\usepackage{cleveref}
\usepackage{mathtools}

\usepackage{authblk}

\usepackage[a4paper,margin=1.1in]{geometry}

\usepackage{amsthm}

\newtheorem{theorem}{Theorem}[section]

\newtheorem{lemma}[theorem]{Lemma}

\theoremstyle{definition}
\newtheorem{definition}[theorem]{Definition}

\theoremstyle{remark}

\usepackage{booktabs} 

\NewDocumentCommand{\pred}{o}{y(t_i, \IfNoValueTF{#1}{\theta}{#1})}
\NewDocumentCommand{\preds}{o}{y\IfNoValueTF{#1}{}{#1}}
\NewDocumentCommand{\real}{}{\mathbb{R}}
\NewDocumentCommand{\mean}{o}{\IfNoValueTF{#1}{\mathbb{E}}{\mathbb{E} \big[#1\big]}} 
\NewDocumentCommand{\var}{o}{\IfNoValueTF{#1}{\operatorname{Var}}{\operatorname{Var} \big(#1\big)}}

\NewDocumentCommand{\data}{m}{\hat{y}_{#1}}
\DeclareMathOperator{\logl}{log\mathcal{L}} 
\NewDocumentCommand{\disc}{o}{(\pred[\IfValueTF{#1}{#1}{\theta}]- \data{i})}
\NewDocumentCommand{\mle}{}{\hat{\theta}_\text{MLE}}

\DeclareMathOperator*{\argmin}{arg\,min}
\DeclareMathOperator*{\tr}{tr}

\title{The Role of Bifurcations in Parameter Estimation: A UQ Analysis of the Non-spatial Klausmeier Model}

\author[1]{Lisa Beer\thanks{\href{mailto:lisa.beer@tum.de}{\texttt{lisa.beer@tum.de}}}}
\author[1]{Christian Kuehn}
\author[1,2]{Chiara Piazzola}

\affil[1]{School of Computation, Information and Technology, Technical University of Munich, Germany}
\affil[2]{ETS de Ingeniería de Caminos, Canales y Puertos, Universitat Politècnica de Catalunya, Barcelona, Spain}

\date{}

\begin{document}

\maketitle
\begin{abstract}
    \input{text/abstract}
\end{abstract}

\vspace{10pt}

\section{Introduction}
\input{text/introduction}

\section{The Klausmeier Vegetation Model}\label{ch:km}
\input{text/klausmeier}

\section{Applied Workflow}\label{sec:workflow}
\input{text/experiment}

\section{Mathematical Tools}\label{sec:tools}
\input{text/methods}

\section{Numerical Experiments}\label{sec:numexp}
\input{text/results}

\section{Conclusion}\label{sec:concl}
\input{text/conclusion}

\newpage
\appendix
    \section{Appendix}
\input{text/appendix}


\printbibliography

\end{document}

%% file: text/abstract.tex
We employ uncertainty quantification methods to investigate how parameter identifiability changes in the vicinity of a bifurcation point.
We perform numerical experiments on the non-spatial Klausmeier vegetation model with random coefficients, which describes biomass-water interactions and exhibits a fold bifurcation.
We partition the domain around the bifurcation into regions with distinct convergence behaviors.
In each region, we follow a UQ workflow that includes sensitivity analysis, Bayesian inference, and Fisher information evaluation to assess parameter identifiability.
The results show that proximity to the bifurcation point is decisive for parameter estimation. Reliable joint inference of both model parameters is not possible when the system exhibits bistability. 
Furthermore, we can reconstruct the model's bifurcation pattern by Fisher information heatmaps, 
presenting a practical tool for bifurcation detection.
Lastly, we highlight the importance of transient data for successful parameter estimation.

%% file: text/introduction.tex
Mathematical models are an indispensable tool in the applied sciences. Analyzing a system's behavior and making predictions about its future course necessitates robust mathematical groundwork, and the underlying tools are often dynamical systems, for instance, in the Life and Earth Sciences:
Dynamical systems can be employed to predict the course of climate phenomena \cite{ClimateModel}, forecast infectious disease outbreaks \cite{SIR}, or model the spread of vegetation patterns \cite{Klausmeier1999}.
To account for parametric variability, it is a standard approach in uncertainty quantification (UQ) to consider dynamical systems within a probabilistic framework and treat the parameters as random variables \cite{UQ2025}. Once the parameter distributions are known, one can proceed with a forward UQ analysis of the dynamics~\cite{breden2020computing,kuehn2021uncertainty,lux2022assessing,hering2026fr}.

The question remains of how to determine reliable parameter distributions.
Parameter reliability is of highest importance to ensure meaningful predictions, particularly around bifurcations~\cite{guckenheimer2013nonlinear}.
Indeed, bifurcations appear under parameter variation and when the parameters cross a critical threshold, the system exhibits different (i.e., topologically non-equivalent) dynamics.
Hence, small uncertainties in parameters near a bifurcation dramatically affect the model's predicted output, highlighting the need for reliable parameter estimates, especially near bifurcations.

This motivates the research in this paper, where we analyze the quality of estimated parameter distributions as a function of their location relative to a bifurcation point.
We perform a case study on a probabilistic version of the non-spatial Klausmeier vegetation model \cite{Klausmeier1999}, which exhibits a fold bifurcation, meaning that two equilibria collide and disappear \cite[Ch. 3.1]{Strogatz2015}.
In numerical experiments, we apply UQ methods \cite{Sullivan2015,UQ2025,UQHandbook} to estimate parameter distributions and analyze their accuracy: We conduct a global sensitivity analysis via Sobol indices \cite{Sobol2001}, and approximate Bayesian posterior distributions \cite{BergerBayes} with Gaussians \cite{InformTheoryFish}. We introduce parameter identifiability and assess it by evaluating parameter correlations and Fisher information \cite{ident25}.

Of particular relevance for this paper is the work by Piazzola et al. \cite{Chiara},
who present a comprehensive overview of UQ methods in epidemiological systems.
Their paper serves as a primary reference, alongside the work by Roesch and Stumpf \cite{Roesch2019}, who investigate parameter inferability around bifurcations in a non-probabilistic setting.
The scope of this paper is to apply the tools of \cite{Chiara} to extend the findings of \cite{Roesch2019}.

The paper is structured as follows: In \cref{ch:km}, we introduce the Klausmeier vegetation model and analyze its bifurcation behavior.
In \cref{sec:workflow}, we describe the workflow and structure of the numerical experiments.
In \cref{sec:tools}, we introduce the mathematical methodology needed. 
In \cref{sec:numexp}, we present our results.
Finally, in \cref{sec:concl} we give a summary and draw a conclusion.

%% file: text/klausmeier.tex
The Klausmeier model \cite{Klausmeier1999} is a PDE model used to describe patterned vegetation spread in semiarid regions like parts of Africa, Australia, and Mexico:
While plants on the flat ground grow in irregular patches, striped patterns are observed if the same species grows on a hillside. The patterns evolve parallel to the hill's contour lines and slowly move uphill over time  \cite{Klausmeier1999,KlausmeierPhd}.
The Klausmeier model consists of the following equations, where $n$ denotes the biomass and  $w$ the downhill flow of surface water \cite{Klausmeier1999}:
\begin{align*}
    \frac{\partial w}{\partial t}&= a -w-wn^2+ v\frac{\partial w}{\partial x}\\
    \frac{\partial n}{\partial t}&= w n^2 - mn + \biggl(\frac{\partial^2}{\partial x^2} + \frac{\partial^2}{\partial y^2}\biggr)n
\end{align*}
The parameter $a>0$ denotes the water input in the form of precipitation, and $m>0$ is the plant mortality. The term $wn^2$ models the water uptake by the vegetation compartment, and $-w$ depicts the baseline loss of water (e.g., evaporation).

For simplicity, we reduce the model to its ODE version. Hence, the previously discussed patterned vegetation states are not further considered, and we focus on scenarios in which the system is in a desert or uniform vegetation state. This reduction provides important insights into the bifurcation structure of the reaction terms. In this study, we thus consider the following system:
\begin{equation}\label{def:klausMod}
    \begin{aligned}
        w'&= a -w-wn^2, \quad w(0) = w_0 \\ 
        n'&= w n^2 - mn, \quad n(0) = n_0
\end{aligned}
\end{equation}

\subsection{Bifurcation Analysis}\label{subs:eqpt}
We are interested in the system's convergence behavior in dependence of the model parameter values.
Hence, we first compute the equilibria of \cref{def:klausMod} to understand where the equations converge to after passing the transient phase \cite[Ch. 1.1]{Wiggins2003}.
To compute the equilibrium states, we set $w'=n'=0$ and solve the following:
\begin{align}
    0 &= a -w-wn^2, \label{app:eq:st1}\\ 
    0 &= w n^2 - mn = n(wn-m) \label{app:eq:st2}
\end{align}
We immediately see that $(w_*, n_*) = (a,0)$ is an equilibrium.
Furthermore, \cref{app:eq:st2} gives the following relation:
\begin{equation}\label{eq:w_eq}
    w = \frac{m}{n}
\end{equation}
Plugging this into \cref{app:eq:st1},
we get $0= a-\frac{m}{n} - mn \Leftrightarrow 0 = mn^2 -an + m $.
Application of the quadratic formula and plugging the result into \cref{eq:w_eq} yields two further equilibria $(w_+, n_+)$ and $(w_-,n_-)$ with
\begin{align}
    w_\pm &= \frac{a\pm\sqrt{a^2-4m^2}}{2}, \label{eq:w+-}\\
    n_\pm &= \frac{a \pm \sqrt{a^2-4m^2}}{2m} \label{eq:n+-}.
\end{align}

\begin{figure}[]
    \centering
    \begin{subfigure}[b]{0.45\textwidth}
        \centering
        \begin{overpic}[width=\textwidth]{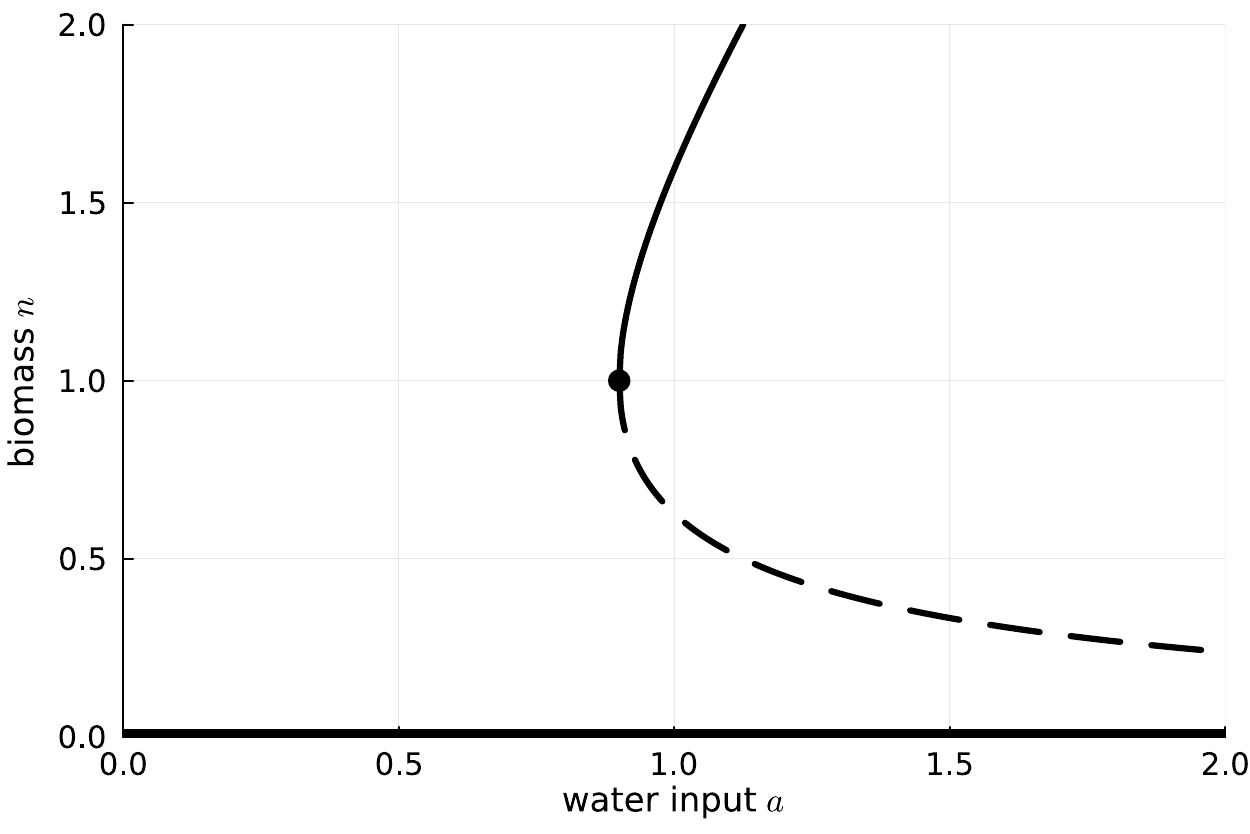}
            \put(44, 0){\color{white}\rule{46pt}{8pt}}
            \put(43, 1){\smaller[2] water input $a$}
            \put(0, 26){\color{white}\rule{8pt}{40pt}}
            \put(0, 26){\rotatebox{90}{\smaller[2] biomass $n$}}
        \end{overpic}
    \end{subfigure} 
    \hfill
    \begin{subfigure}[b]{0.45\textwidth}
        \centering
        \begin{overpic}[width=\textwidth]{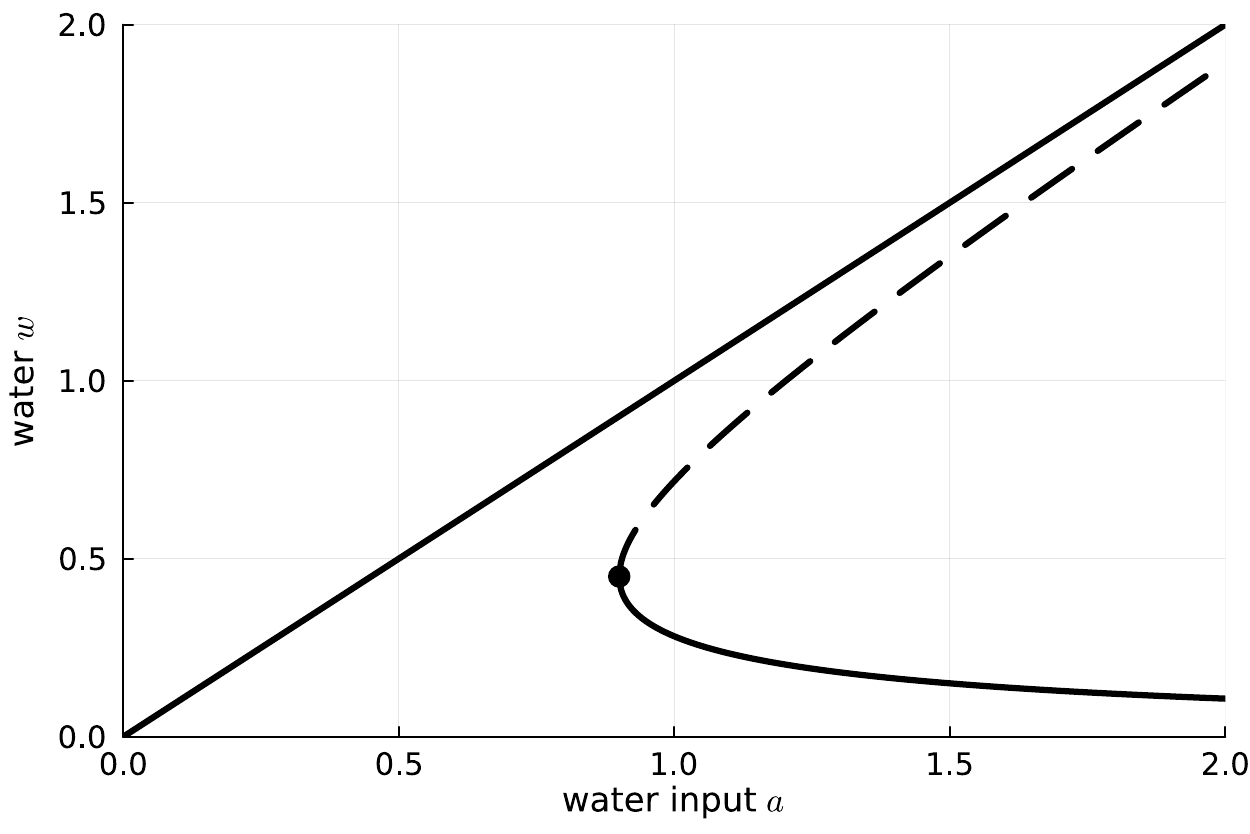}
            \put(44, 0){\color{white}\rule{46pt}{8pt}}
            \put(43, 0){\smaller[2] water input $a$}
            \put(0, 26){\color{white}\rule{8pt}{40pt}}
            \put(0, 28){\rotatebox{90}{\smaller[2] water $w$}}
        \end{overpic}
    \end{subfigure}
    \caption{Bifurcation diagrams of $w$ and $n$ for $m=0.45$ with free bifurcation parameter $a$.}
    \label{fig:bifurcation}
\end{figure}
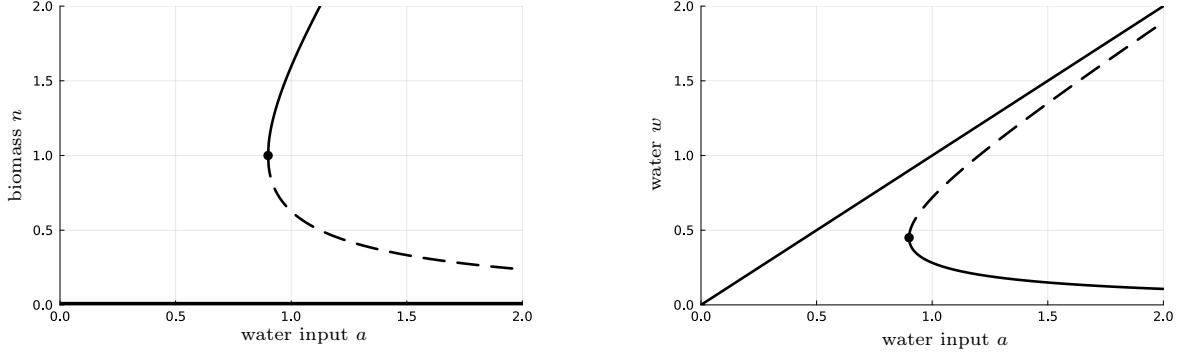

If we analyze \cref{def:klausMod}, we observe that \cref{eq:w+-,eq:n+-} possess no real solutions for $a<2m$ and two distinct, real solutions for $a>2m$.
Therefore, two distinct equilibria collide and disappear at $a=2m$, which is defined as a fold bifurcation \cite[Ch. 3.1]{Strogatz2015}. Such fold bifurcations occur along the parameter subspace
\[\{
(a,m)^\top \in \mathbb{R}^2 : 2m = a
\}.\]
We visualize the equilibria and bifurcation points in an exemplary 2D-bifurcation diagram  in \cref{fig:bifurcation} \cite[Ch. 3.1]{Strogatz2015}. We use $m=0.45$ as a parameter value for the plant mortality \cite{Klausmeier1999}, while $a$ serves as free bifurcation parameter.

For $a>2m$ we have two distinct equilibrium curves. One of these is stable, i.e., it attracts nearby trajectories, while the other one is unstable, i.e., it repels nearby trajectories.
The solid lines in \cref{fig:bifurcation} represent the stable equilibrium, the dashed lines the unstable one\footnote{The stability of the equilibria is derived in a linear stability analysis, see \cref{sub:stabana}.}.
We see that for $a>2m$, the system exhibits multistability (specifically, bistability), with two distinct attracting steady states.

\paragraph{Interpretation of the bifurcation diagram:}
In the bifurcation diagrams we consider, we focus on the steady-states (or equilibrium points). Of particular practical importance are asymptotically stable steady states, where a system settles after passing an initial transient phase.
For instance, if we consider the point $(a^*,w^*)^\top = (1,1)^\top$ on the solid equilibrium line in \cref{fig:bifurcation}, $w^*$ is the respective steady state value of the $w$ compartment for a model simulation with parameter $a^*$ and a suitable initial condition $w_0$. Suitable means that $w_0$ lies within the basin of attraction of the respective equilibrium branch.
In this specific example, a suitable $w_0$ must lie above the solid line or in between the solid line and the dashed line, see \cref{fig:bifurcation}, as the dashed line represents a repelling state in the phase space, which turns out in our case to be a border between different basins of attraction for asymptotically stable steady states.

In \cref{sec:div_prm_sp}, we make use of this interpretation of the bifurcation diagram to divide the parameter and initial condition space around the bifurcation into regions exhibiting distinct convergence behaviors.

\subsection{Introducing Parametric Uncertainty}\label{sec:rand_km}
We introduce uncertainty into the Klausmeier vegetation model by treating the parameters $a$ and $m$ and the initial conditions $ w_0$ and $ n_0$ as random variables.
With this approach, we account for the variability that model parameters express in the real world \cite{UQHandbook,UQsmith,UQluis}.

On the probability space $(\Omega, \mathcal{F}, \mathcal{P})$, the Klausmeier model with random parameters becomes
\begin{equation}\label{def:herz}
    y'= f\big(t, y, \theta(\omega)\big), 
    \quad y\big(t_0,\theta(\omega)\big) = y_0(\omega), \quad  \omega\in\Omega
\end{equation}
where we define $y'\coloneq \big(w',n'\big)^\top$, $y_0(\omega)\coloneq \big(w_0(\omega),n_0(\omega)\big)$, $\theta(\omega)\coloneq\big(a(\omega),m(\omega)\big)^\top$ and $f(t,y,\theta(\omega))\coloneq\big(\theta_1(\omega) - w - w n^2, w n^2 - \theta_2(\omega) n\big)^\top$.
Sometimes we write $\theta$ instead of $\theta(\omega)$ for brevity.

%% file: text/experiment.tex
This section motivates and explains the numerical experiments we conduct in \cref{sec:numexp} and anticipates the mathematical tools needed.
We consider the random Klausmeier model in \cref{def:herz} as the basis for the analysis.
The workflow is organized in steps.

\subsection{Step 0: Division of Parameter Space}\label{sec:div_prm_sp}
The final goal of this work is to assess how proximity to a bifurcation point affects our ability to estimate parameters.
Therefore, we first classify a parameter's location with respect to the bifurcation by suitably subdividing the parameter and initial condition space.

\begin{figure}[]
    \centering
    \begin{subfigure}[b]{0.45\textwidth}
        \centering
        \begin{overpic}[width=\textwidth]{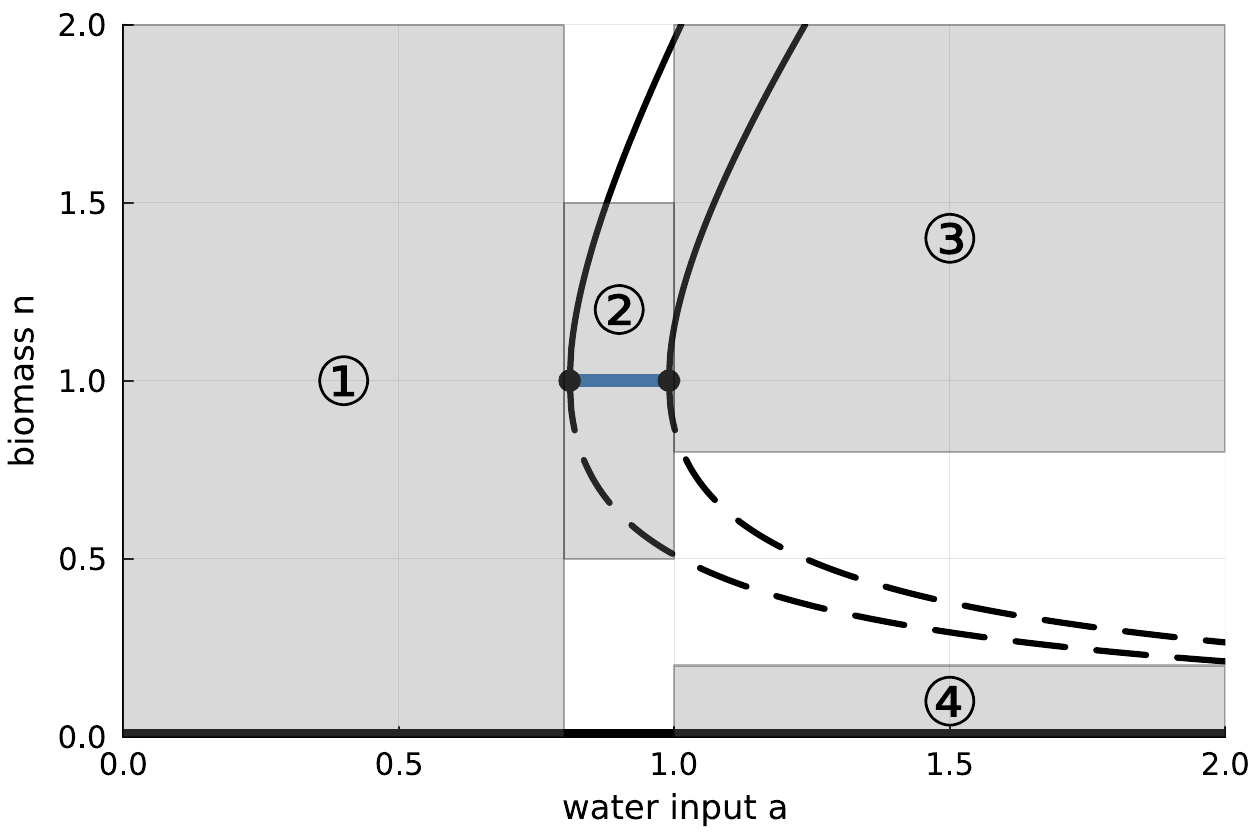}
            \put(44, 0){\color{white}\rule{46pt}{8pt}}
            \put(43, 1){\smaller[2] water input $a$}
            \put(0, 26){\color{white}\rule{8pt}{40pt}}
            \put(0, 34){\rotatebox{90}{\smaller[2] $n_0$}}
        \end{overpic}
    \end{subfigure} 
    \hfill
    \begin{subfigure}[b]{0.45\textwidth}
        \centering
        \begin{overpic}[width=\textwidth]{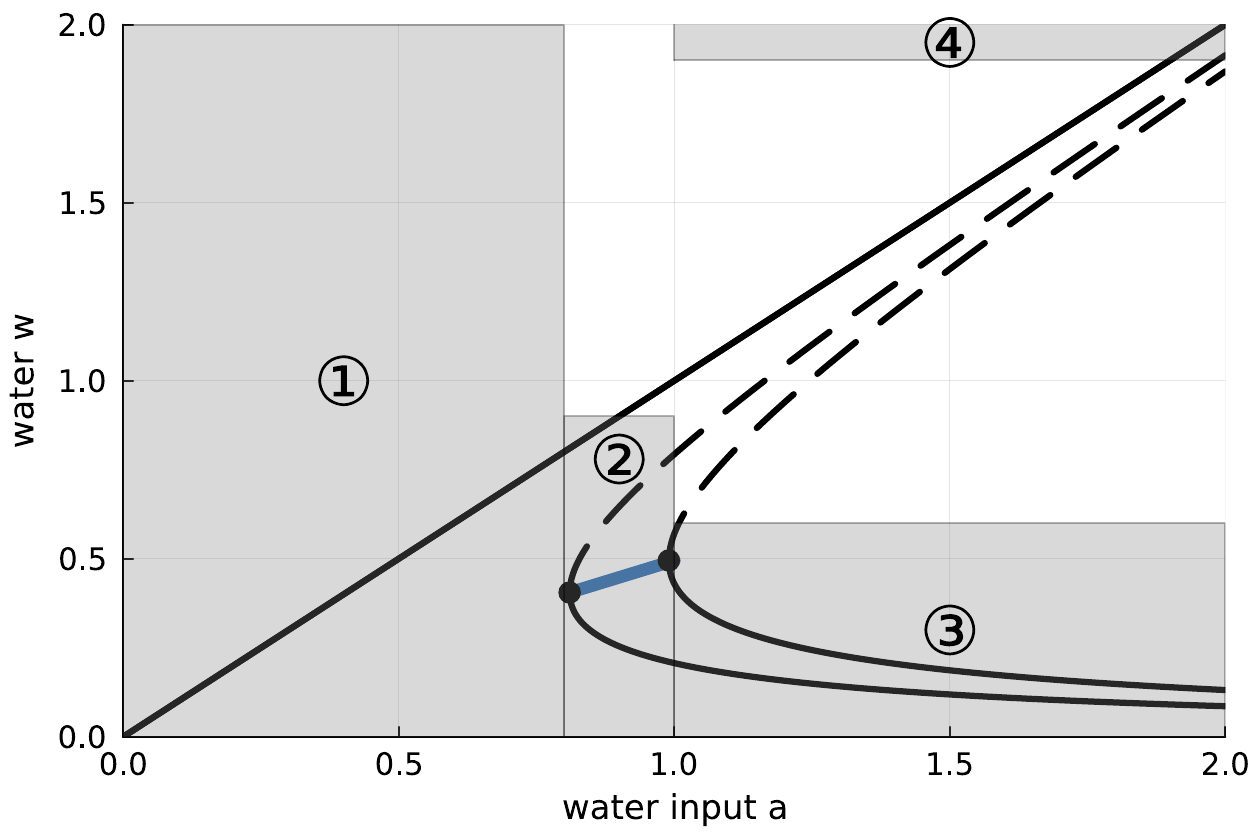}
            \put(44, 0){\color{white}\rule{46pt}{8pt}}
            \put(43, 0){\smaller[2] water input $a$}
            \put(0, 26){\color{white}\rule{8pt}{40pt}}
            \put(0, 34){\rotatebox{90}{\smaller[2] $w_0$}}
        \end{overpic}
    \end{subfigure}
    \caption[]{Classifying regions around the bifurcation.
    Each region represents part of a different basin of attraction.
    The bifurcation point varies with $m$, as indicated by the blue line.}
    \label{fig:prm_regions}
\end{figure}
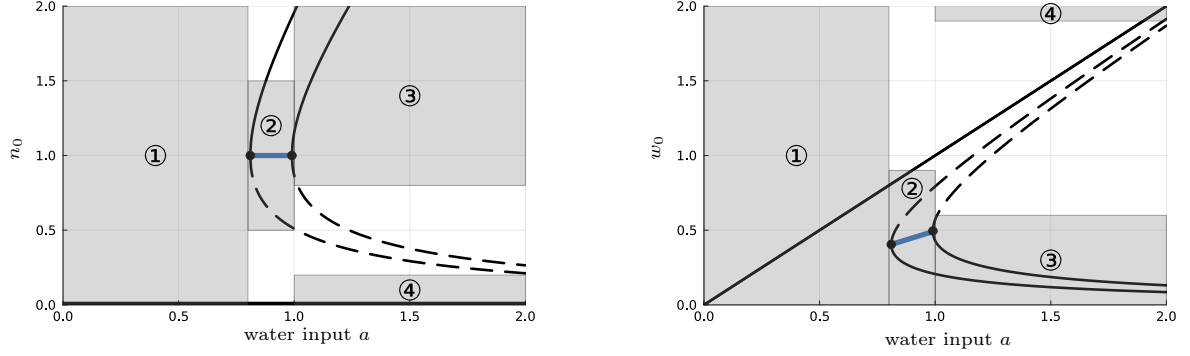

We differentiate four regions, as indicated in \cref{fig:prm_regions}.
The grey boxes describe the potential values of parameter $a$, and the initial conditions $w_0, n_0$.
Furthermore, we allow $m$ to vary in $[0.405,0.495]$.
We overlay the grey boxes in \cref{fig:prm_regions} by the equilibrium branches resulting from choosing $m$ as its extremal values (for the intermediate values $m\in(0.405,0.495)$ the bifurcation point lies on the blue line and the equilibrium branches vary accordingly).

With this plot, it becomes visible that the grey boxes each represent different regions of attraction in the parameter space, according to the interpretation of the bifurcation diagram in \cref{subs:eqpt}:
The system converges to the plant-free equilibrium for all parameter and initial condition combinations from regions one and four, whereas it converges to the nontrivial equilibrium for all parameter and initial condition combinations from region three. In region two, we observe an alternating convergence behavior for different parameter and initial conditions combinations, as the location of the bifurcation and hence the attractor strongly depends on the parameter value of $m$.

\subsection{Step 1: Sobol Sensitivity Analysis}\label{sec:sob_wf}
Before addressing the parameter estimation task in step 2, we investigate the model's sensitivity to parameter variations in each region.
In doing so, we hope to gain valuable information about the parameters, since sensitivity results are closely linked to parameter inferability: if a parameter causes little sensitivity in the model output, it is probably difficult to estimate it using data-driven methods. 

Throughout the sensitivity analysis, we consider the random Klausmeier model introduced in \cref{sec:rand_km} and assume that the parameters and initial conditions follow uniform distributions corresponding to the parameter regions introduced in \cref{sec:div_prm_sp}.
In all subspaces shown in \cref{fig:prm_regions}, we analyze the time-dependent output variance as parameters vary. This reveals, in which region which parameters or initial conditions cause the most sensitivity.
Specifically, we conduct a global sensitivity analysis using Sobol indices. For more information on Sobol indices, see \cref{sec:sobol}.
The experiments and their results are presented in \cref{sec:sob_res}. 

\subsection{Step 2: Bayesian Inverse Analysis}\label{sec:inv_wf}
After the sensitivity analysis, we address the inverse problem.
For each region $r, r\in\{1,\dots,4\}$, we pick a representative parameter point $\theta^r = (a^r,m^r)^\top$ together with some initial conditions $y_0^r=(w_0,m_0)^\top$ and simulate trajectory data adding mean-zero Gaussian noise.
We try to estimate the joint distribution of $a$ and $m$ via Bayesian inversion, see \cref{sec:invuq}.
We approximate the posterior with a Gaussian centered at the MLE, see \cref{sec:gaus_app}.
Finally, we determine the accuracy of the MLE and Gaussian approximation at each parameter point $\theta^r$ by investigating its identifiability, see \cref{sec:ident}. 
We obtain a local result of parameter inferability for each region $r$ at the respective point and initial condition combination $\theta^r, y_0^r$.

We perform two experiments. In the first one, we observe data from the transient phase and the steady state of the system. In the second, we assume that no transient-phase data are available. By comparing these two experiments, we assess the role of transient phase information in parameter estimation. The experiments and their results are presented in \cref{sec:pi_res}.

\subsection{Step 3: Identifiability Landscapes via Fisher Information}\label{sec:glob_fish}
We are interested not only in inferability within the distinct regions but also in a more global picture of parameter inferability around the bifurcation.
Therefore, we choose the Fisher information of a parameter estimate as the primary measure of identifiability, see \cref{sec:fi}. 
We consider a point grid of $a, n_0$ combinations for fixed $m$ and $w_0$.
For each point, we simulate data, add mean-zero Gaussian noise, and compute the MLE. 
Similar to \cite{Roesch2019}, we then create heatmaps by mapping the Fisher information of the MLE to the underlying point $(a,n_0)^\top$.
This approach reveals regions where the Fisher information of the MLEs is high, indicating greater identifiability.
We repeat the experiment at different noise levels and with varying numbers of simulated data points.
The experiments and their results are presented in \cref{sec:fi_res}.

%% file: text/methods.tex
This section introduces the necessary mathematical methods. We first discuss sensitivity analysis, followed by inverse UQ methods.

Without loss of generality, we consider some one-dimensional model $y'=f\big(t,y,\theta(\omega)\big)$ with an $m-$ dimensional random parameter vector $\theta$. On a discretized time scale, we characterize the predicted model output at $t_i$, $i\in \{0,\dots,M\}$, given the parameter vector $\theta(\omega)\in\real^m$ by the tuple $\big(t_i,\pred\big)$. 
Furthermore, we assume $\theta \sim U[0,1]^m$ throughout the section.

\subsection{Global Sensitivity Analysis}\label{sec:sobol}
With sensitivity analysis, we aim to link variability in the model's output to uncertainty in its input parameters, i.e., we want to quantify the influence that variation in a parameter within its distribution has on the model output \cite{Saltelli2008}.

We conduct a global sensitivity analysis that reveals information across the entire parameter input space.
To be precise, we use Sobol indices as a sensitivity measure.

\begin{definition}[First-order Sobol index \protect{\cite[p.161]{Saltelli2008}}]\label{def:sob-index}
    The first-order Sobol index of a parameter $\theta_k$ at time $t_i$ is 
    \begin{equation*}
        S_k(t_i) \coloneq \frac{\var[E\big[\preds(t_i, \theta)\big|\theta_k\big]]}{\var[\pred]}.
    \end{equation*}
\end{definition}
The formula computes the variance of the expected model output, provided the parameter $\theta_k$ is fixed. Dividing by the total variance yields the proportion of the output variance that arises purely from the parameter $\theta_k$. Hence, the larger $S_k(t_i)$, the greater the influence of $\theta_k$ on the model output. We see that $0\leq S_k(t_i)\leq 1$ \cite[Ch. 4.4]{Saltelli2008}.

Furthermore, we are interested in the total Sobol indices, which consider not only the parameter $\theta_k$ itself but also its interactions with all other parameters.
\begin{definition}[Total Sobol index \protect{\cite[p. 163]{Saltelli2008}}]\label{def:TO-sob-index}
    The total Sobol sensitivity index of $\theta_k$ is defined as
    \begin{equation*} 
        S_{T_k}(t_i) \coloneq 1 - \frac{\var[E\big[\pred\big| \theta_{\sim k}\big]]}{\var[\pred]}, 
    \end{equation*}
    where $\theta_{\sim k} \coloneq (\theta_1, \dots, \theta_{k-1},\theta_{k+1},\dots,\theta_m)$.
\end{definition}
The formula subtracts the variance of the expected model output, with all parameters except $\theta_k$ held fixed, from the total output variance. Hence, it measures exactly all the variance arising from $\theta_k$ and its interaction with other parameters.
By definition it follows that $S_k(t_i) \leq S_{T_k}(t_i) \leq 1$ for all $k$ \cite[Ch. 4.5]{Saltelli2008}.

If we evaluate the Sobol indices not only for a fixed time point $t_i$, but instead consider several time points $t_0 < t_1 < \cdots < T$ on an interval $[0,T]$, we are able to express the model sensitivity as time-dependent trajectories.

The Sobol indices can be computed using numerical Monte Carlo-based methods, see \cite[4.6]{Saltelli2008}.

\subsection{Inverse Uncertainty Quantification}\label{sec:invuq}
In this section, we address how to determine parameter distributions that fit the observed model output data.
Therefore, we assume that we have synthetic data $\hat{y}_i$ at time points $t_i$ for $i\in \{0,\dots,M\}$ stemming from a true underlying parameter realization $\theta^*$.
The data is corrupted by additive Gaussian noise and is modeled by
\begin{equation}\label{eq:syn_dat}
    \data i = \pred[\theta^*] + \varepsilon_i, \quad \varepsilon_i\sim  \mathcal{N}(0, \sigma^2)\hspace{4pt} \text{i.i.d}.
\end{equation}
The noisiness makes the standard inverse problem of finding a single true parameter estimate for $\theta^*$ ill-posed, motivating once more the need to employ probabilistic and UQ methods \cite[Ch. 1.1]{UQluis}\cite[Ch. 6.1]{Sullivan2015}.

Our goal is to estimate a suitable distribution for $\theta^*$.
We address this problem with Bayesian inversion \cite{BergerBayes}.
Thereby, we first assume some parameter distribution, called "prior" distribution, which is then updated to a "posterior" distribution taking into account the information obtained by the synthetic data.
The posterior distribution describes which parameter values are most plausible given the data \cite{BergerBayes}.
This is formalized in Bayes' theorem:
\begin{theorem}[Bayes' theorem \protect{\cite[Ch. 4]{BergerBayes}\cite[Ch. 8.1]{UQsmith}}]\label{thm:bayes}
    Let $\theta$ follow some prior distribution $\pi(\theta)$, 
    let $(\data1,\dots, \data M)^\top$ be noisy data.
    We obtain a data-informed update of $\pi(\theta)$ by
    \begin{equation}\label{eq:posterior}
        \pi(\theta \mid \data1,\dots, \data M) = \frac{p(\data1,\dots,\data M| \theta)\pi(\theta)}{\displaystyle\int p(\data1,\dots,\data M| \theta)\pi(\theta)d\theta},
    \end{equation}
    where $p(\data1,\dots, \data M| \theta)$ is the data distribution.
    We call $\pi(\theta\mid\data1,\dots,\data M)$ the posterior distribution. It describes how probable different values of $\theta$ are conditioned on the observed data.
\end{theorem}
If we view the joint probability density function of all data observations as a function of the parameter $\theta$, we obtain the likelihood function $\mathcal{L}(\theta)\coloneq p(\data1,\dots, \data M| \theta).$ 
Applying the logarithm does not qualitatively change the likelihood surface, but simplifies computation.
The log-likelihood function in our specific setup is then \cite{Stortelder1998}:
\begin{equation*}
    \logl(\theta) \propto -\frac{1}{2 \sigma^2} \sum_{i=1}^{M}(\pred - \data i).
\end{equation*}
The (log-)likelihood function quantifies how well the model predictions under a $\theta$ agree with the observation data.
In the considered case (i.i.d. Gaussian noise and uniform priors), the best possible choice for $\theta$ is the maximum likelihood estimate (MLE) \cite{BergerBayes,Held2021Likeli}:
\begin{equation*}
    \mle \coloneq \argmin_{\theta \in [0,1]^m} \biggl(-\logl(\theta)\biggl).
\end{equation*}

\subsubsection{Gaussian Approximation}\label{sec:gaus_app}
Approximating the posterior by a Gaussian can simplify its form and facilitate, e.g., sampling.
In our setting, such an approximation is possible locally around the MLE, if $-\nabla^2\logl(\hat{\theta}_\text{MLE})$ is invertible, and the log-likelihood function is three times continuously differentiable.
The approximation is then \cite[Result 8]{BergerBayes}:
\begin{equation}\label{thm:gaus_appr}
        \hat\pi(\theta\mid \data1,\dots, \data M) = \mathcal{N}\big(\hat\theta_{\mathrm{MLE}},(\nabla^2\logl(\hat\theta_{\mathrm{MLE}}))^{-1}\big)
\end{equation}
To see that this holds, we rewrite the posterior distribution on the parameter domain to 
\begin{equation}\label{eq:post_in_gaus_pf}
    \begin{aligned}
    \pi(\theta \mid \data1, \dots, \data M) 
    &= \frac{\exp\big(\logl(\theta)\big)}{\displaystyle \int \exp\big(\logl(\theta)\big)\, d\theta} \, \mathbf{1}_\Theta(\theta),
    \end{aligned}
\end{equation}
where we use \cref{eq:posterior} and the prior $\theta\sim U[0,1]^m$.
We now compute the second order Taylor expansion of $\logl(\theta)$ at the point $\hat{\theta}_\text{MLE}$:
    \begin{equation}\label{eq:taylor}
        \logl(\theta) = \logl(\mle) + \frac{1}{2}(\theta-\mle)^\top  \left.\nabla^2\logl(\theta)\right|_{\theta = \mle} (\theta-\mle).
    \end{equation}
If we now plug \cref{eq:taylor}
into \cref{eq:post_in_gaus_pf}, we obtain the local approximation
$\hat\pi(\theta\mid\data1,\dots,\data M)\approx\pi(\theta\mid\data1,\dots,\data M)$ as
\begin{equation*}
    \hat\pi(\theta\mid\data1,\dots,\data M) = 
    \frac{\exp\biggl(\frac{1}{2}(\theta-\mle)^\top\left.\nabla^2\logl(\theta)\right|_{\theta = \mle} (\theta-\mle)\biggr)}
    {(2\pi)^{m/2}\det\big(\left.\nabla^2\logl(\theta)\right|_{\theta = \mle}^{-1}\big)^{1/2}},
\end{equation*}
$\forall\theta\in\Theta$.
Clearly, $\hat\pi(\theta\mid\data1,\dots,\data M)$ is the multivariate Gaussian from \cref{thm:gaus_appr}.

The Gaussian approximation is a local approximation and not always appropriate in practice, as it crucially depends on the MLE and the covariance matrix. In \cref{sec:fi,sec:ident}, we will discuss the cause and the impact of an ill-shaped mode of the Gaussian approximation and how this relates to reliable parameter estimation.

\subsubsection{Fisher Information}\label{sec:fi}
The covariance matrix of the Gaussian approximation (\ref{thm:gaus_appr}) is exactly the inverse Fisher information matrix $I(\theta)$,
\begin{equation*}
    I(\theta) \coloneq -\nabla^2\logl(\theta).
\end{equation*}
The scalar Fisher information is the trace $\tr I(\theta)$.
By definition, via the Hessian, the Fisher information measures the log-likelihood's local curvature around its evaluation point. Hence, $\tr I(\mle)$ describes the accuracy of the MLE:
A peaked log-likelihood surface, i.e., a high Fisher information value, indicates that the parameter estimate is very reliable, as a slight change in the estimate causes a significant decrease in the log-likelihood and thus a less plausible parameter estimate.
Similarly, a flat log-likelihood surface, i.e., a low Fisher information value, indicates that the maximum likelihood estimate is uncertain, as the surrounding estimates create a similar log-likelihood output, and, thus, serve as equally plausible parameter estimates \cite{Roesch2019}.

The connection between the Fisher information matrix and the covariance matrix illustrates that low Fisher information values lead to inaccuracies in the Gaussian approximation, as they indicate high variance. This can lead to the loss of identifiability, as discussed in the next section.

\subsubsection{Identifiability}\label{sec:ident}
We already anticipated in \cref{sec:fi,sec:gaus_app} that the validity of the MLE and the posterior's Gaussian approximation can be compromised in practice.

Indeed, on the one hand, high parameter correlations can lead to ill-conditioned covariance matrices for the approximating Gaussian posteriors.
On the other hand, high correlations are reflected in stretched log-likelihood level sets, where the change of one parameter can be compensated for by another.
Both imply that parameters cannot be estimated independently of each other \cite[Sec. 4]{Miao}.

Since parameter correlations are computed by normalizing the covariance matrix $C \coloneq I(\mle)^{-1}$, we see that these considerations are linked to the Fisher information.
Hence, both the evaluation of the Fisher information and the assessment of parameter correlations provide valuable insights into the quality of the estimated posteriors.

All these considerations are summarized under the umbrella of identifiability.
In a noiseless, ideal setting, we say that a parameter $\theta^*$ is structurally identifiable if it is uniquely determined in an inverse analysis of the model output:
\begin{equation*}
    y(t,\theta) = y(t,\theta^*)\quad\forall t\in[0,T] \quad\Longrightarrow\quad \theta=\theta^*.
\end{equation*}
Structural identifiability is necessary for successful parameter estimation, yet it is not sufficient.
As discussed in the previous paragraphs, practical identifiability is not guaranteed and may be undermined by noisy and limited data \cite[Sec. 7., 8.]{Chiara}\cite{wielandStrident,ident25}.

%% file: text/results.tex
In this section, we present and discuss the experiments in the workflow from \cref{sec:workflow}.
The underlying code is written in the Julia programming language.
Model simulations are done with the \texttt{DifferentialEquations.jl} \cite{DiffEq_jl} package.
To implement the computation of Sobol indices, we use \texttt{GlobalSensitivityAnalysis.jl} \cite{GlobalSensitivityAnalysis_jl}.
For the computation of MLEs, the package \texttt{Optim.jl} \cite{Optim_jl} is used.
To obtain the Hessian, e.g. needed for the computation of the Fisher information, we use the package \texttt{ForwardDiff.jl} \cite{ForwardDiff_jl}.
We present the most interesting results below.
The complete experiments and the code are available on \href{https://github.com/lisab00/ma-code}{GitHub}.

\subsection{Step 1: Sobol Sensitivity Analysis}\label{sec:sob_res}
We first perform the Sobol analysis motivated in \cref{sec:sob_wf}.
We uniformly draw $5000$ parameter samples from each of the four regions introduced in \cref{sec:div_prm_sp} and compute solution trajectories for $t_i\in[0,100]$ for each sample. We then compute the first-order and total Sobol indices (see \cref{sec:sobol}) at $1000$ equidistant time points to obtain time-dependent sensitivity trajectories.

\begin{figure}
    \centering
    \begin{subfigure}[b]{0.45\textwidth}
        \centering
        \begin{overpic}[width=\textwidth]{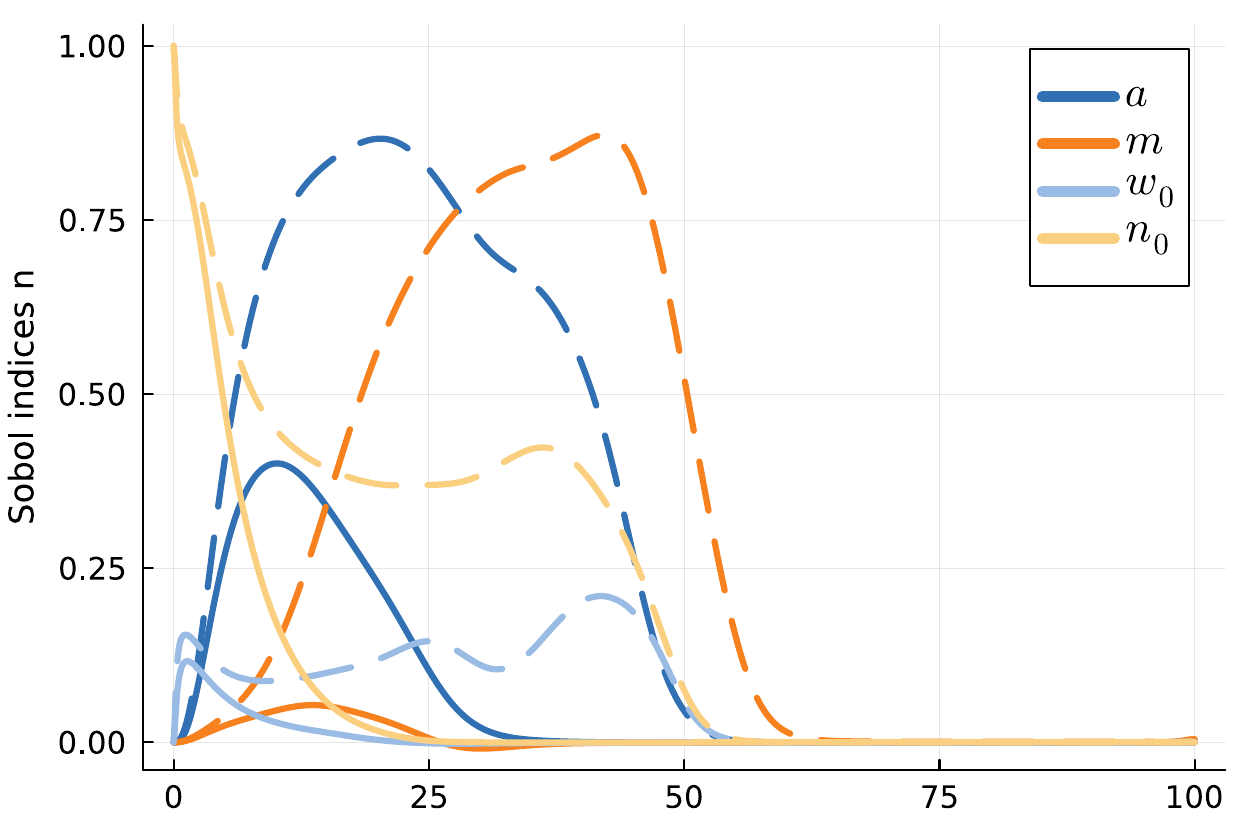}
            \put(48, -3){\color{white}\rule{40pt}{8pt}}
            \put(49, -1){\smaller[2] time}
            \put(0, 20){\color{white}\rule{8pt}{80pt}}
            \put(-1, 21){\rotatebox{90}{\smaller[2] Sobol indices $n$}}
        \end{overpic}
    \end{subfigure}
    \hfill
    \begin{subfigure}[b]{0.45\textwidth}
        \centering
        \begin{overpic}[width=\textwidth]{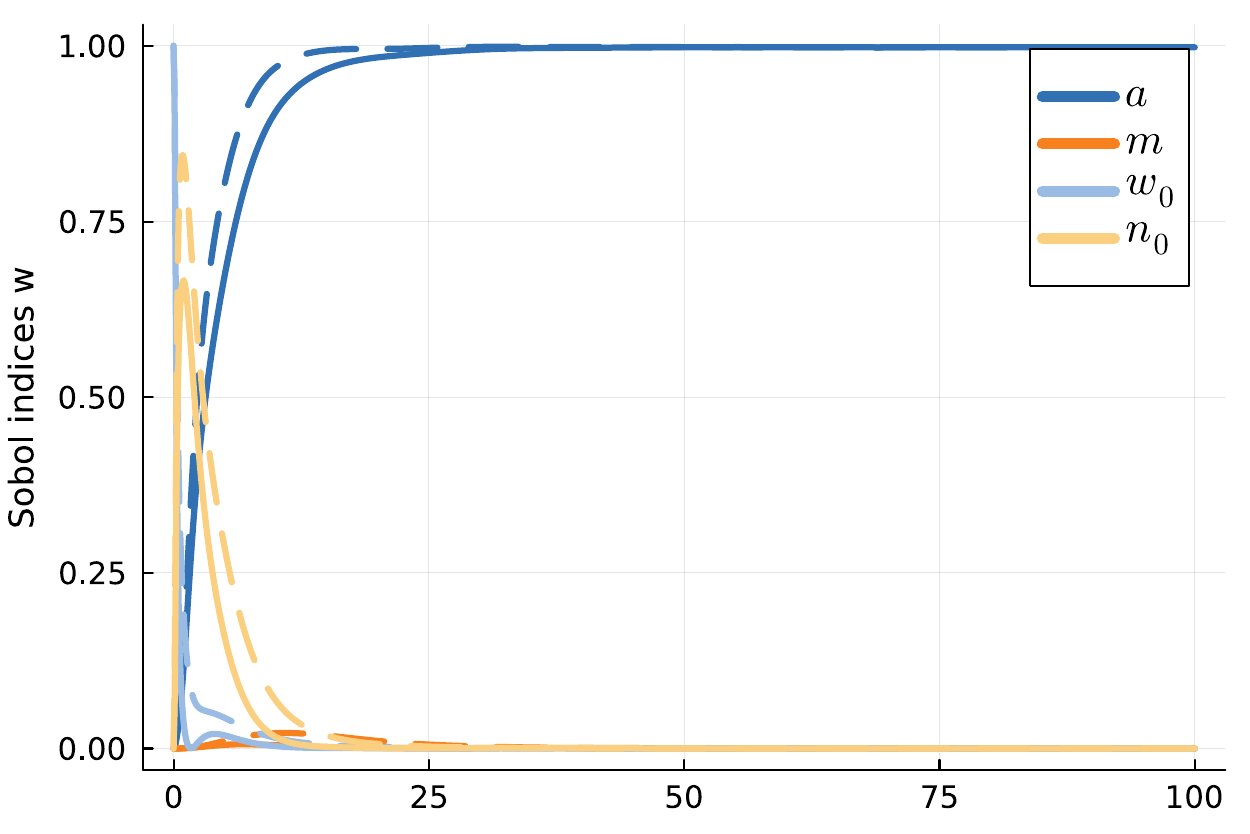}
            \put(48, -3){\color{white}\rule{40pt}{8pt}}
            \put(49, -3){\smaller[2] time}
            \put(0, 20){\color{white}\rule{8pt}{80pt}}
            \put(-1, 21){\rotatebox{90}{\smaller[2] Sobol indices $w$}}
        \end{overpic}
    \end{subfigure}
    \caption[Time-dependent Sobol indices for priors in parameter region one.]{Time-dependent Sobol indices for priors in parameter region one. Solid lines: first-order index, dashed lines: total indices.
}
    \label{fig:region1_sobol}
\end{figure}
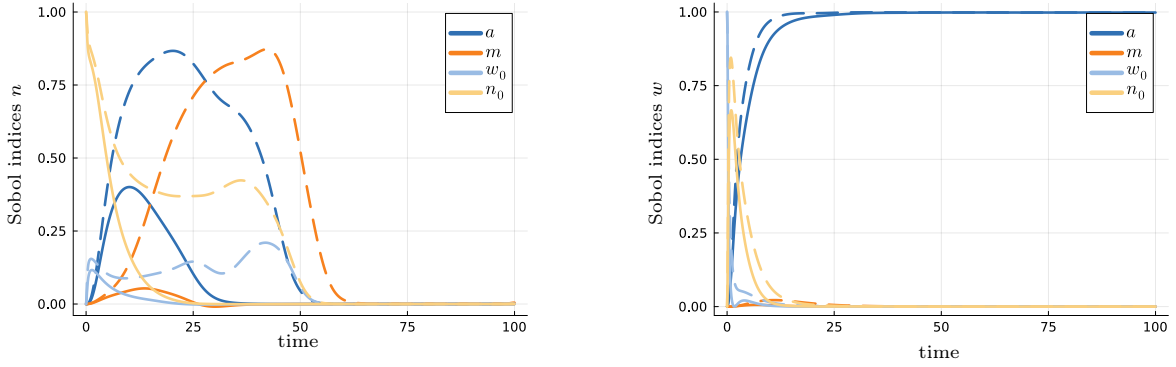

\begin{figure}
    \centering
    \begin{subfigure}[b]{0.45\textwidth}
        \centering
        \begin{overpic}[width=\textwidth]{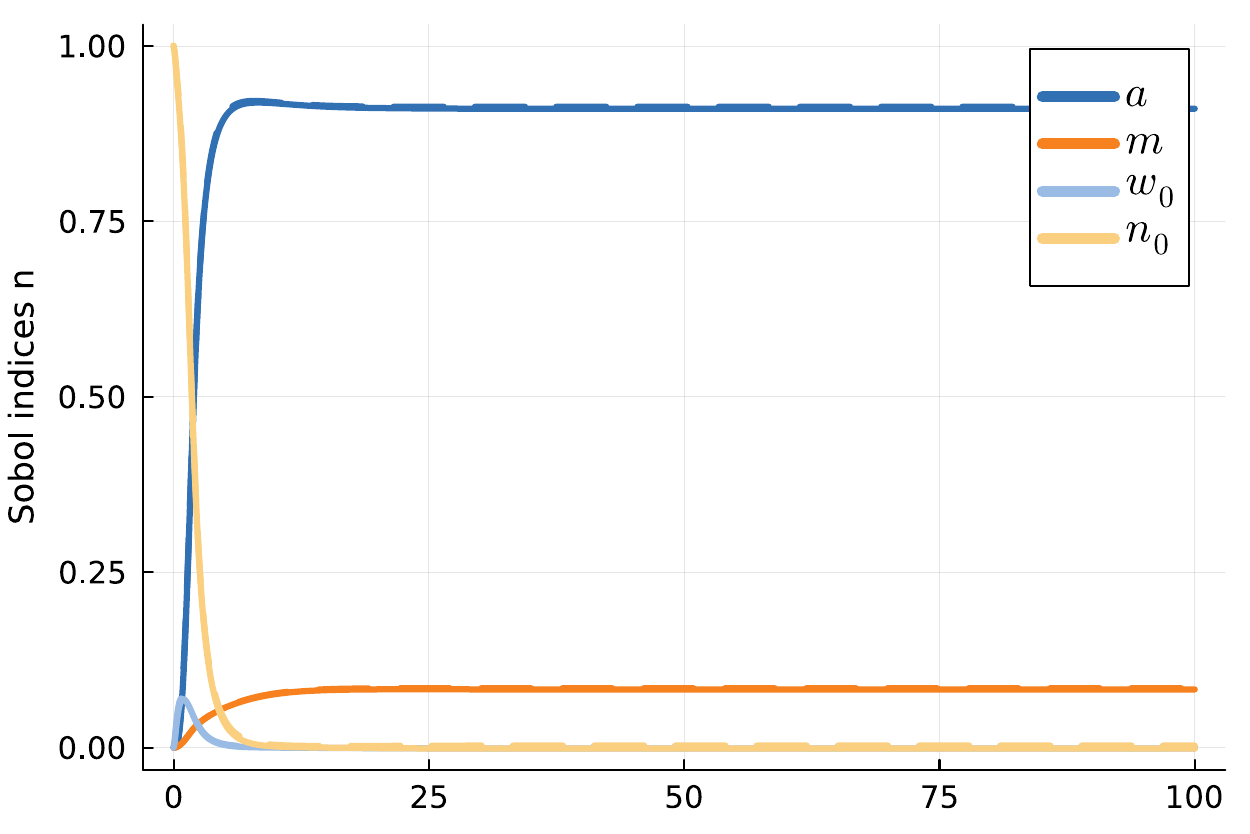}
            \put(48, -3){\color{white}\rule{40pt}{8pt}}
            \put(49, -3){\smaller[2] time}
            \put(0, 20){\color{white}\rule{8pt}{80pt}}
            \put(-1, 21){\rotatebox{90}{\smaller[2] Sobol indices $n$}}
        \end{overpic}
    \end{subfigure}
    \hfill
    \begin{subfigure}[b]{0.45\textwidth}
        \centering
        \begin{overpic}[width=\textwidth]{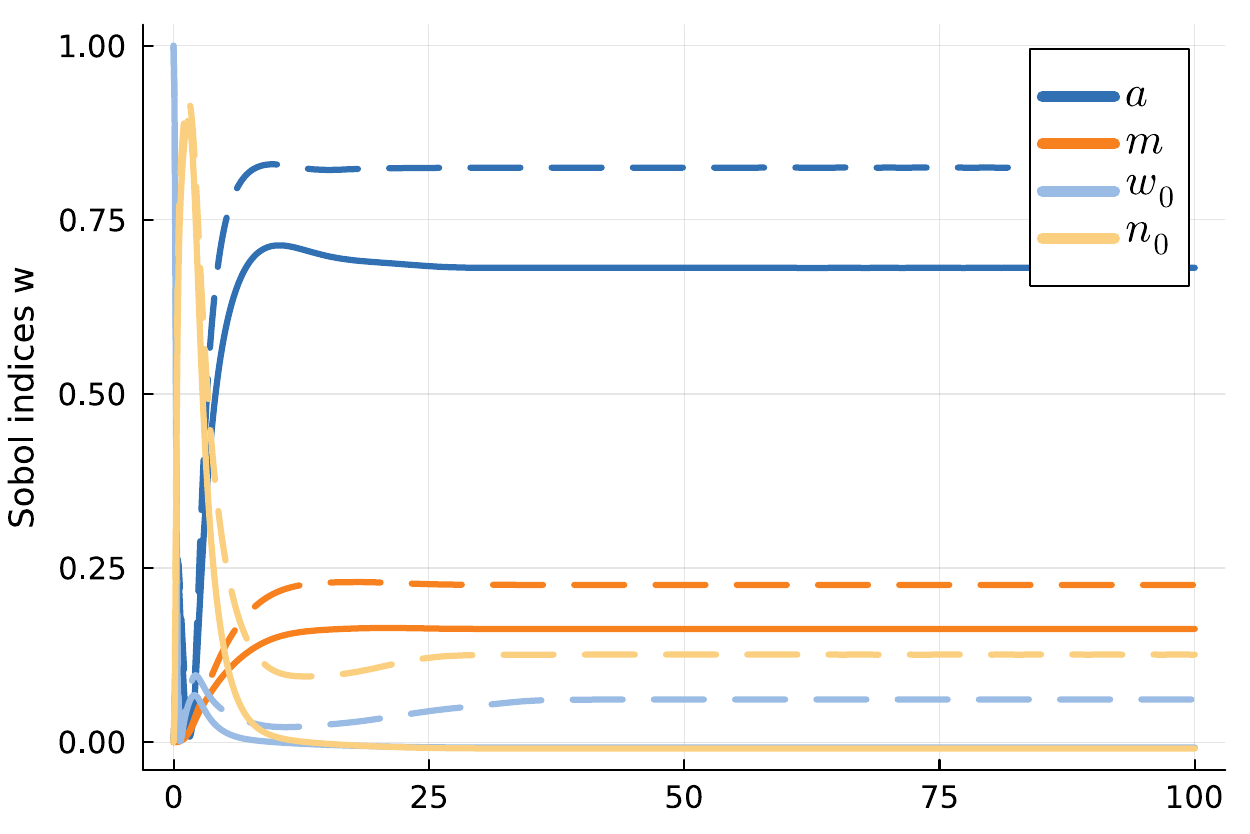}
            \put(48, -3){\color{white}\rule{40pt}{8pt}}
            \put(49, -3){\smaller[2] time}
            \put(0, 20){\color{white}\rule{8pt}{80pt}}
            \put(-1, 21){\rotatebox{90}{\smaller[2] Sobol indices $w$}}
        \end{overpic}
    \end{subfigure}
    \caption[Time-dependent Sobol indices for priors in parameter region three.]{Time-dependent Sobol indices for priors in parameter region three. Solid lines: first-order index, dashed lines: total indices.
}
    \label{fig:region3_sobol}
\end{figure}
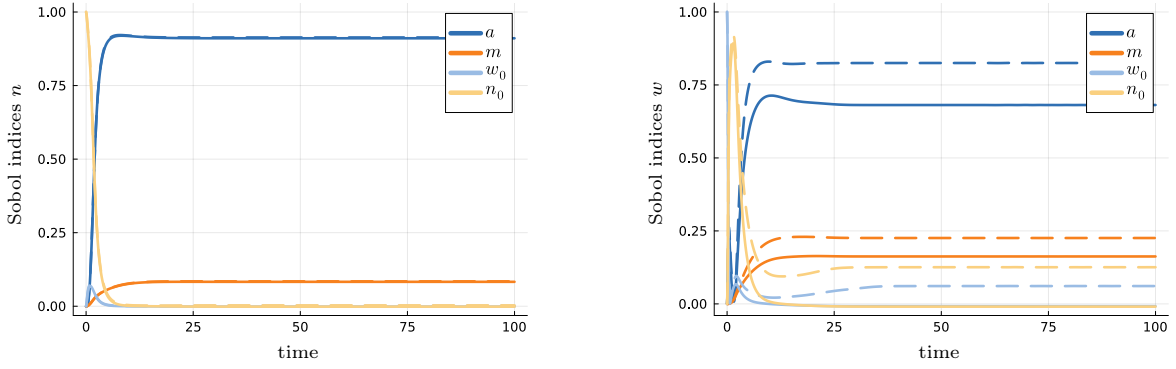

We focus here on the results for regions one and three, see \cref{fig:region1_sobol,fig:region3_sobol}.
The behavior in region four is similar to region one, and in region two, the conclusions are less significant, as the alternating convergence behavior is naturally reflected in the model output sensitivity.

In region one, see \cref{fig:region1_sobol}, the system converges to the plant-free equilibrium.
Once in the stable state, the value of the $w$ compartment is solely determined by $a$ (we see in \cref{fig:bifurcation,fig:prm_regions} that the linear relationship $w=a$ holds), which is why all sensitivity of $w$ concentrates in $a$, and higher order interactions play no significant role.
In the $n$ compartment, all parameter and initial condition values lead to the same stable state value, namely zero.
Hence, the sample variance converges rapidly to zero, making the Sobol indices technically undefined. We report them as zero here, consistent with the fact that variation in parameters and initial conditions causes no model output sensitivity. We see that in the transient phase, particularly parameter interactions cause high sensitivity, before all indices settle to zero.

In region three, see \cref{fig:region3_sobol}, the system converges to the non-trivial equilibrium.
From \cref{fig:bifurcation,fig:prm_regions}, it becomes clear that the steady-state values for both compartments depend heavily on $a$ and $m$, as reflected in the sensitivity caused by these parameters. We note that $a$ causes a significantly higher sensitivity than $m$.
Moreover, it is interesting to note that in the $w$ compartment, interactions of $w_0,n_0$ also play a significant role, although their first-order indices are negligible.

In both regions one and three, we detect that the variation of particularly $a$, but also $m$, causes sensitivity.
Hence, it seems reasonable that we can estimate both, or at least one, in step 2.
Meanwhile, the initial conditions $w_0$ and $n_0$ appear to cause little to no sensitivity in the model output, especially once the stable state is reached. Hence, we keep them fixed at some value for the numerical experiments in steps 2 and 3, unless stated otherwise.

Finally, we note that the duration of the transient phase is reflected in the time-dependent Sobol indices, which converge once the system has settled into a steady state.

\subsection{Step 2: Bayesian Inverse Analysis}\label{sec:pi_res}
We now address step 2 as explained in \cref{sec:inv_wf}.
Based on the results from the preceding Sobol analysis, we fix the initial conditions $w_0$, $n_0$, and hope to be able to jointly estimate both parameters $a$ and $m$.

The synthetic data is generated using both compartments at $M=100$ equidistant time steps, with mean-zero Gaussian noise of variance $\sigma^2=0.1$.

\subsubsection{Experiment 1: Including Transient Phase Information}
First, we conduct the experiments for the case where the system is observed during the transient phase and in its stable state ($t_0=0, t_\text{end}=100$).

\begin{figure}[]
    \centering
    \begin{subfigure}[b]{0.45\textwidth}
        \begin{overpic}[width=\textwidth]{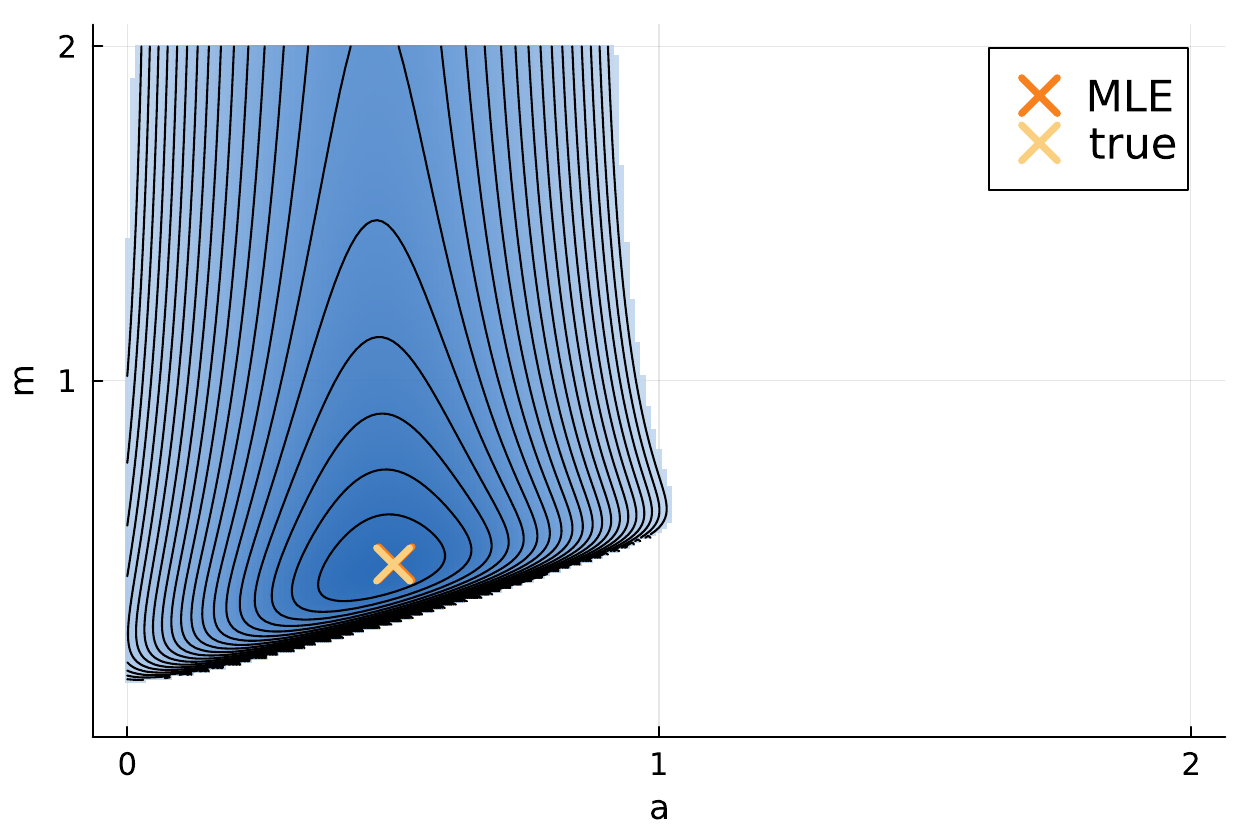}
            \put(48, 0){\color{white}\rule{40pt}{8pt}}
            \put(50, 0){\smaller[2] $a$}
            \put(0, 20){\color{white}\rule{8pt}{80pt}}
            \put(0, 33){\rotatebox{90}{\smaller[2] $m$}}
        \end{overpic}
        \caption{Log-likelihood surface.
    Values below a certain threshold are mapped to white to improve visibility.}
    \label{fig:pt1_ll}
    \end{subfigure}
    \hfill
    \begin{subfigure}[b]{0.45\textwidth}
        \begin{overpic}[width=\textwidth]{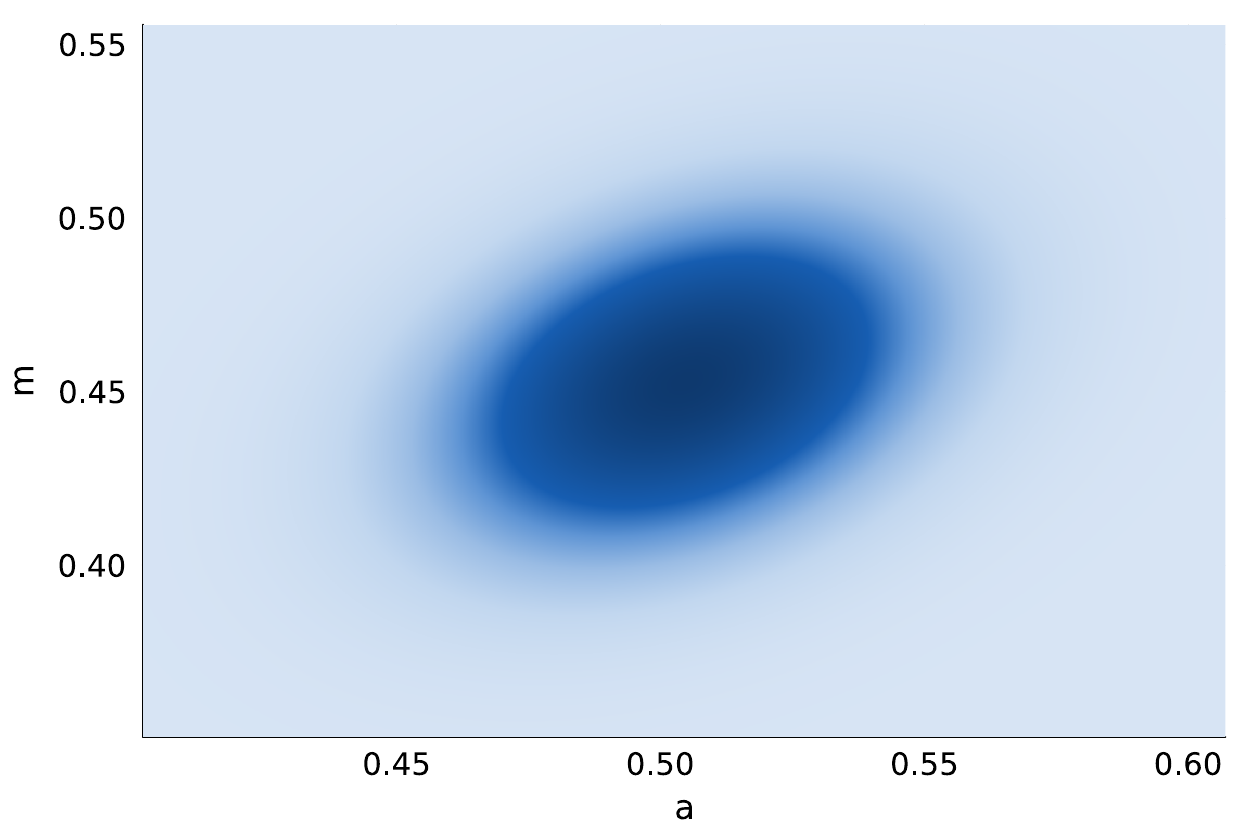}
            \put(48, 0){\color{white}\rule{40pt}{8pt}}
            \put(50, 0){\smaller[2] $a$}
            \put(0, 20){\color{white}\rule{8pt}{80pt}}
            \put(0, 33){\rotatebox{90}{\smaller[2] $m$}}
        \end{overpic}
        \caption{Contour plot of the Gaussian density used as posterior approximation.}
        \label{fig:pt1_gh}
    \end{subfigure}
    \caption{Region one: Results for $\theta_1$, estimation of $a$, $m$.}
    \label{fig:pt1_combined}
\end{figure}

\Cref{fig:pt1_ll} shows the log-likelihood surface and the obtained MLE in region one. We see that the contours are well-shaped and possess a clear maximum. The MLE is indistinguishable from the true  parameter.
Analyzing the correlation of the parameter estimates, we conclude that both $a$ and $m$ are identifiable in this case and the Gaussian approximation represents their posterior density adequately (for the exact numerical values, we refer to the results on \href{https://github.com/lisab00/ma-code/tree/main/notebooks/identifiability/identifiability_of_points}{GitHub}). A visualization of the approximating Gaussian density can be found in \cref{fig:pt1_gh}. The shape is barely stretched in one direction, highlighting the small parameter correlation and equal magnitudes of the variances of both covariates.
\begin{figure}[]
    \centering
    \begin{subfigure}[b]{0.45\textwidth}
        \begin{overpic}[width=\textwidth]{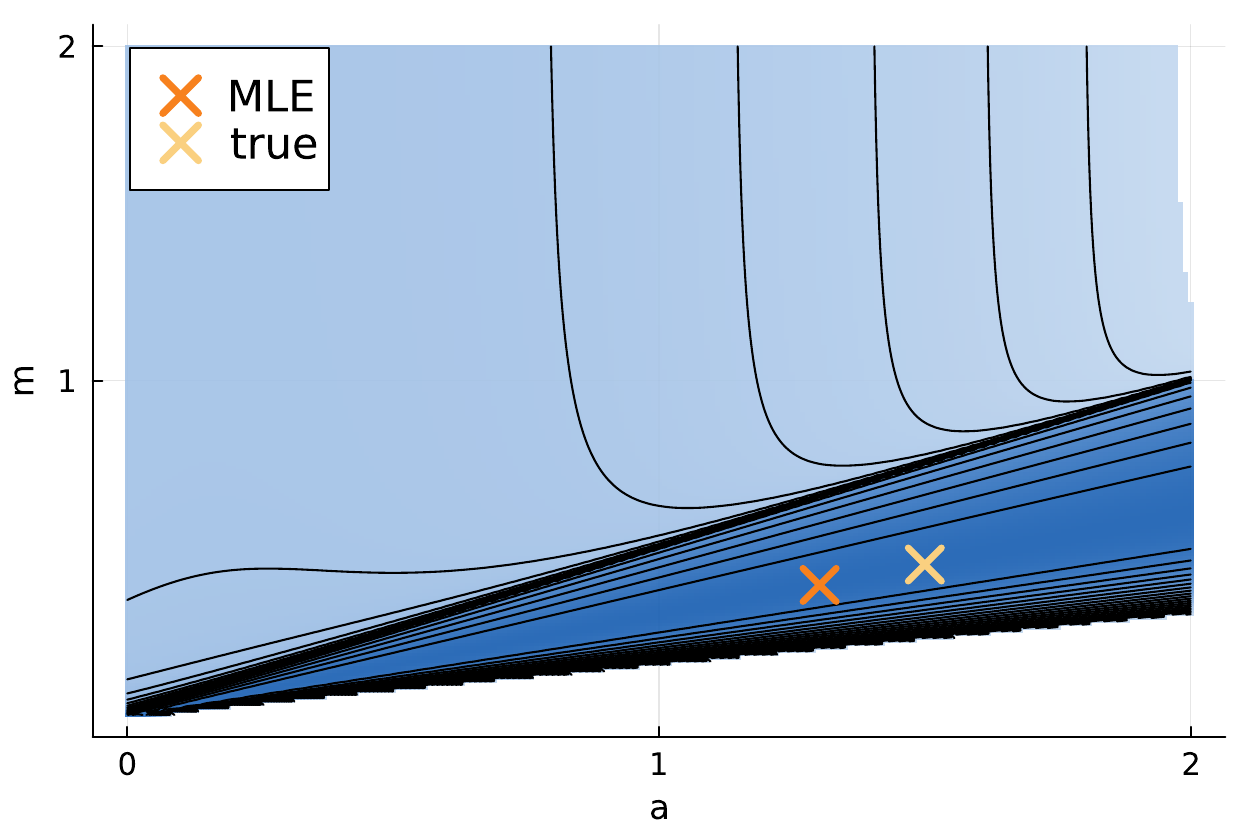}
        \put(48, 0){\color{white}\rule{40pt}{8pt}}
        \put(50, 0){\smaller[2] $a$}
        \put(0, 20){\color{white}\rule{8pt}{80pt}}
        \put(0, 33){\rotatebox{90}{\smaller[2] $m$}}
    \end{overpic}
    \caption{Log-likelihood surface. Values below a certain threshold are mapped to white to improve visibility.}
    \label{fig:pt3_ll}
    \end{subfigure}
    \hfill
    \begin{subfigure}[b]{0.45\textwidth}
        \begin{overpic}[width=\textwidth]{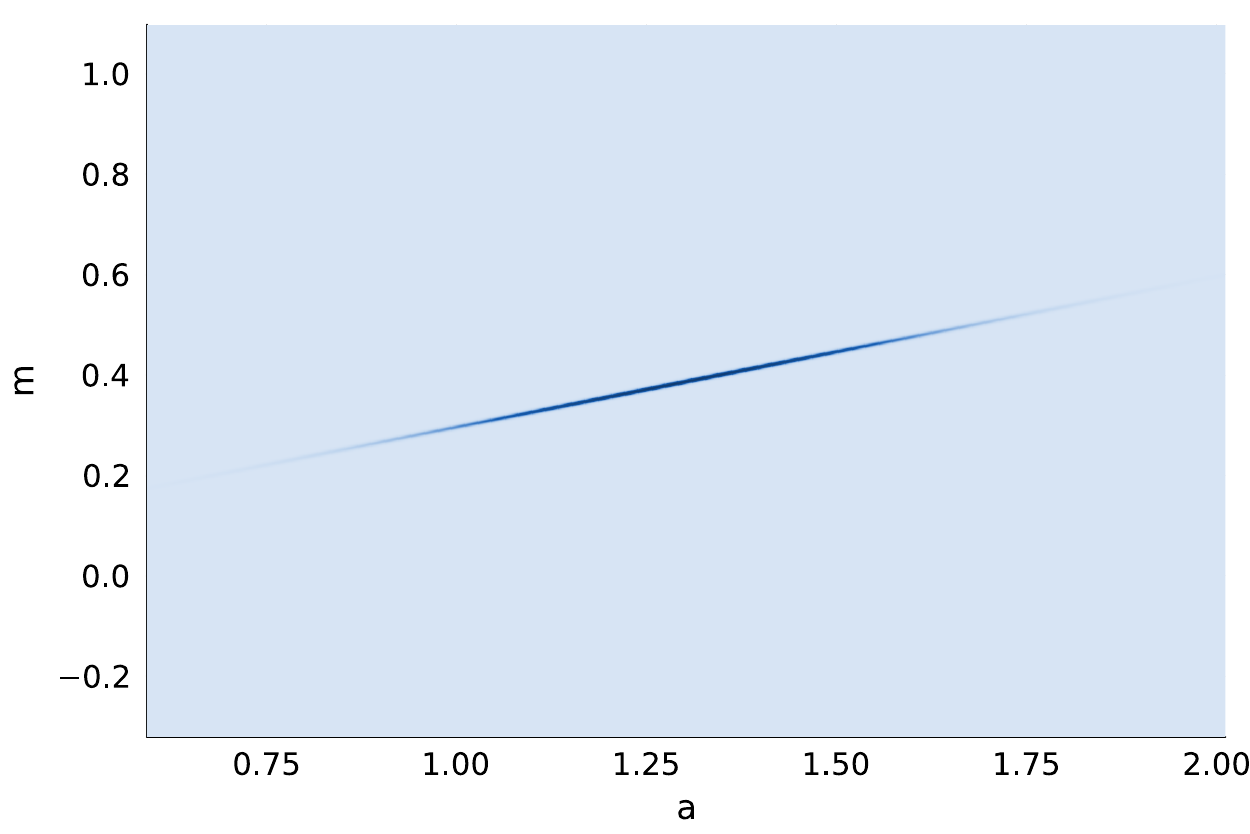}
        \put(48, 0){\color{white}\rule{40pt}{8pt}}
        \put(50, 0){\smaller[2] $a$}
        \put(0, 20){\color{white}\rule{8pt}{80pt}}
        \put(0, 33){\rotatebox{90}{\smaller[2] $m$}}
    \end{overpic}
    \caption{Contour plot of the Gaussian density used as posterior approximation.}
    \label{fig:pt3_gh}
    \end{subfigure}
    \label{fig:pt3_combined}
    \caption[]{Region three: Results for $\theta_3$, estimation of $a$, $m$.}
\end{figure}
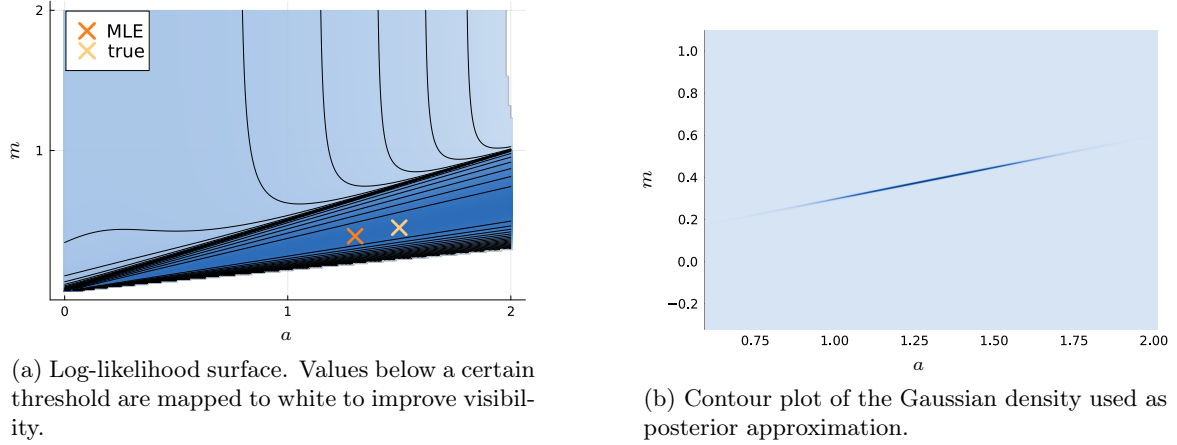

Further experiments reveal that joint estimation of the parameters $a$ and $m$ is possible only in region one. In the other regions, issues arise due to strong parameter correlations.
As an example, the log-likelihood surface of the parameter estimates for region three is shown in \cref{fig:pt3_ll}.
The contour lines are almost flat in one direction, and a clear maximum cannot be distinguished.
This is also reflected in the approximating Gaussian density, see \cref{fig:pt3_gh}. The scales indicate that the variance of $a$ is significantly higher than that of $m$. The components of the MLE are highly correlated and exhibit near-linear dependence. Hence, the parameters are practically not identifiable, and the Gaussian approximation is not suitable as an approximation of the posterior.

Hence, in contrast to the expectations based on the sensitivity analysis results, both $ a$ and $ m$ are identifiable only in region one. We resort to estimating only $a$, as it exhibits the greatest sensitivity in the Sobol sensitivity analysis.
We find that $a$ is identifiable across all parameter regions, see the experiments on \href{https://github.com/lisab00/ma-code/tree/main/notebooks/identifiability/identifiability_of_points}{GitHub}.

\subsubsection{Experiment 2: Observation of Stable State Only}
In a second experiment, we consider the case where no transient phase data is available and we have only steady-state data ($t_0=100, t_\text{end}=200$).
The estimation of $a$ remains comparable in quality across all regions, but the lack of transient phase information makes the simultaneous inference of $m$ impossible, even in region one.
In particular, this reflects well the conclusions drawn from the Sobol analysis of region one (see \cref{fig:region1_sobol}), where we observed that sensitivity with respect to $m$ is exhibited only within the transient phase.

\subsection{Step 3: Identifiability Landscapes via Fisher information}\label{sec:fi_res}
\begin{figure*}[]
    \centering
    
    \begin{tikzpicture}
    
        \node (myplot) at (0,0) {
            \includegraphics[width=0.6\textwidth]{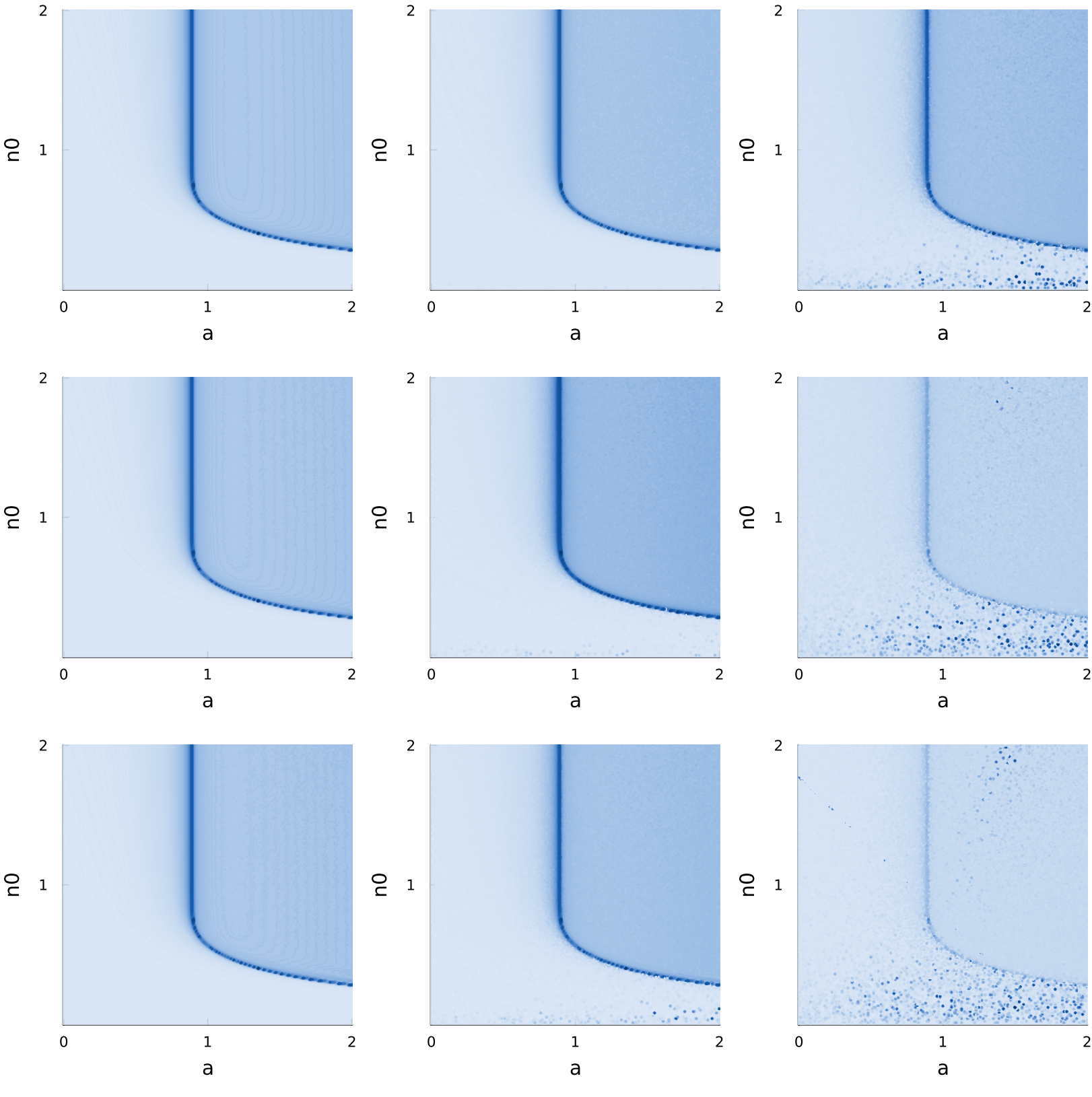}
        };
        
        \draw[thick, -stealth] ($(myplot.north west) + (1.5, 0.1)$) -- ($(myplot.north east) + (-1.5, 0.1)$)
              node[midway, above=0mm, font=\smaller] {increase noise level};

        \draw[thick, -stealth] ($(myplot.north west) + (-0.2, -1.5)$) -- ($(myplot.south west) + (-0.2, 1.5)$)
              node[pos=0.2, left=2mm, rotate=90, font=\smaller] {decrease number of data points};
              
    \end{tikzpicture}
    
    \caption{Visualization of the FI at the MLE for $a$, $m$.
    Dark colors correspond to a high FI.}
    \label{fig:fi_grid_am} 
\end{figure*}
\begin{figure*}[]
    \centering
    
    \begin{tikzpicture}
    
        \node (myplot) at (0,0) {
            \includegraphics[width=0.6\textwidth]{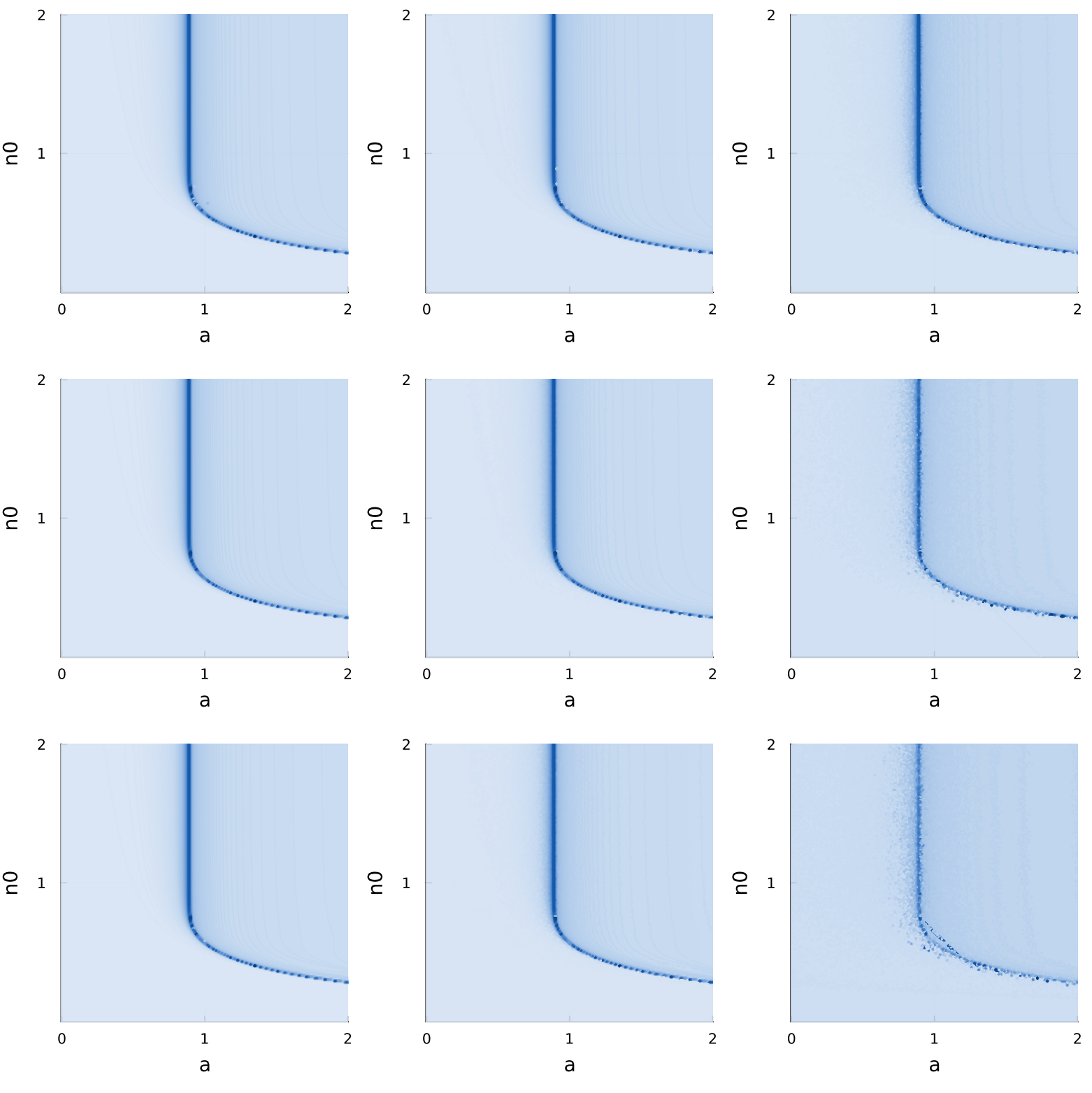}
        };
        
        \draw[thick, -stealth] ($(myplot.north west) + (1.5, 0.1)$) -- ($(myplot.north east) + (-1.5, 0.1)$)
              node[midway, above=0mm, font=\smaller] {increase noise level};

        \draw[thick, -stealth] ($(myplot.north west) + (-0.2, -1.5)$) -- ($(myplot.south west) + (-0.2, 1.5)$)
              node[pos=0.2, left=2mm, rotate=90, font=\smaller] {decrease number of data points};
              
    \end{tikzpicture}
    
    \caption{Visualization of the FI at the MLE for $a$.
    Dark colors correspond to a high FI.}
    \label{fig:fi_grid_a} 
\end{figure*}

We now want to visualize the Fisher information of parameter estimates around the bifurcation as explained in \cref{sec:glob_fish}.
We fix $m=0.45$, $w_0=1.0$, and vary $a,n_0\in[0,2]$ around the bifurcation.
For each evaluation point, we simulate the model trajectories ($t_0 =0, t_\text{end}=100$), add mean-zero Gaussian noise, and infer the joint MLE of $a$ and $m$ (or, only $a$). We vary the number of data points $M\in\{50,100,500\}$ and the noise level $\sigma^2\in\{0.01,0.1,0.5\}$ to simulate real-world conditions. 
We then evaluate the Fisher information at the MLE and visualize it in heatmaps.

The results for the joint estimation of $a,m$ are presented in \cref{fig:fi_grid_am}, and the results for the estimation of $a$ only are presented in \cref{fig:fi_grid_a}.
We observe that values separating the attraction regions yield substantially higher Fisher information, as even minor changes to the parameters $a$ and $m$ (or only $a$) can cause the simulation to converge to a different stable state. Therefore, the likelihood function is highly sensitive to parameter variations, yielding high Fisher information values.
In both cases, we can clearly visualize the model's bifurcation structure.
 
Comparing different noise levels and data amounts, we find that both have little effect on the Fisher information values when estimating only $a$, see \cref{fig:fi_grid_a}. This explains why our estimation results for $a$ are consistently positive across all regions, whereas when estimating $a$ and $m$, we obtained more distinct results from case to case.
In general, we see that the increase in noise has a considerably more severe impact on the Fisher information values than the decrease in sample size.
The effect is most apparent in region four, where the Fisher information values appear more random than structured.
In contrast, region one is the least affected by noise and data decrease.
This explains why inference for both parameters at point $\theta_1$ was successful, see \cref{fig:pt1_combined}, especially for the chosen initial condition $n_0=1.5$.

At first, the high Fisher information values around the bifurcation point in \cref{fig:fi_grid_am} seem to contradict the fact that the joint parameter estimation of $a,m$ in these regions is not possible. However, this phenomenon is an artifact of reducing the Fisher information matrix to its trace. For instance, a high diagonal value corresponding to $a$ can lead to a large trace of the Fisher information matrix, even though the diagonal value corresponding to $m$ is close to zero.
Hence, relying solely on the Fisher information as an identifiability indicator can be misleading, as reducing it to a scalar inevitably leads to a loss of information.

\subsection{Discussion}\label{sec:sum_res}
The numerical experiments highlight that reliable estimation of both model parameters $a$ and $m$ is only possible in parameter region one, provided that transient data are available. In the other regions, practical identifiability fails due to strong parameter correlations. This is reflected in linearly stretched level sets of the log-likelihood and, hence, in an inadequate Gaussian approximation of the posterior, see \cref{sec:pi_res}.

We conclude that identifiability decreases when the data stems from the bistable regime. Therefore, in practice, high parameter correlations could serve as an indicator for multistability regions. Parameter identifiability, hence, not only depends on the attraction region itself, but also on the number of attractors present in the corresponding parameter regime.
This is illustrated when comparing the results for regions one and four. Although both exhibit convergence to the plant-free equilibrium and the Sobol results are similar, the identifiability results differ.

Crucially, the convergence in region four is very rapid. This once more highlights the importance of transient data, as underscored by the two experiments in step 2, \cref{sec:pi_res}.

Furthermore, we could reconstruct the model's bifurcation pattern from Fisher information heatmaps, see \cref{sec:fi_res}. These heatmaps indicate subspaces in which parameter estimation is possible under real-world conditions.
This analysis reveals the power of Fisher information for bifurcation detection via model simulation, but also highlights that combining multiple identifiability measures is essential for a comprehensive picture.

%% file: text/conclusion.tex
In this paper, we investigated how the identifiability of the parameters in the Klausmeier ODE model changes around its fold bifurcation.

We could show that successful parameter estimation heavily relies on the parameter location and the availability of transient data. Furthermore, we were able to reconstruct the model's bifurcation pattern in a Fisher information heatmap.

The key approach of our work was to divide the parameter and initial condition space appropriately by interpreting the bifurcation diagram with respect to different attraction regions.
We concluded that joint estimation of both parameters is not possible in the bistable regime after the bifurcation, and that we would need historical data from a non-vegetated state for parameter calibration.
Together with the importance of transient data, as highlighted by our experiments, this underscores the need for long-term data collection in real-world applications.

A limitation of this work is that the Klausmeier model is very simple, having only two compartments and two parameters to calibrate.
Still, we hope that the fundamental insights obtained here can be transferred to more complex models that exhibit similar, although more complex, bifurcation behavior.

%% file: text/appendix.tex
\subsection{Klausmeier Model Stability Analysis}\label{sub:stabana}
For the stability analysis of \cref{def:klausMod}, we use \cref{thm:stable_ev}:
\begin{theorem}[\protect{\cite[Thm 1.2.5]{Wiggins2003}}]\label{thm:stable_ev}
    Consider a system $y'=f(y)$ with an equilibrium point $y^*$.
    If all eigenvalues $\lambda$ of the Jacobian $Jf(y_*)$ satisfy $\Re(\lambda) < 0$, then $y_*$ is stable.  
    If at least one $\Re(\lambda) > 0$, then $y_*$ is unstable.
\end{theorem}
In the following, assume $a>2m$.
The Jacobian is
\begin{equation*}
    J(w,n) = \biggl(\nabla\frac{\partial w}{\partial t},
        \nabla\frac{\partial n}{\partial t}\biggr)^\top
    = 
    \begin{pmatrix}
        -1-n^2 & -2wn \\
        n^2 & 2wn-m
    \end{pmatrix}.
\end{equation*}
Using \cref{thm:stable_ev} we immediately conclude that the equilibrium $(w_*,n_*)$ is asymptotically stable, as $$J(w_*,n_*) = J(a,0) = \begin{pmatrix}
    -1 & 0 \\
    0 & -m
\end{pmatrix}$$
is in diagonal form and has only real, negative eigenvalues for all $a, m\in\mathbb{R}_+$.

For the stability analysis of the remaining equilibria we exploit that $m = w_\pm n_\pm$ (see \cref{eq:w_eq}) and transform the Jacobian at the equilibria to
\begin{equation*}
    J(w_\pm, n_\pm) = 
    \begin{pmatrix}
        -1-n_\pm^2 & -2m \\
        n_\pm^2 & m
    \end{pmatrix}.
\end{equation*}
Its eigenvalues are: 
\[\lambda_\pm = -\frac{1}{2}(1+n_\pm^2 -m) \pm \sqrt{\frac{(1+n_\pm^2 -m)^2}{4} + m(1-n^2_\pm)}\]
Setting $\alpha \coloneq -\frac{1}{2}(1+n_\pm^2 -m)$ and $\beta \coloneq m(1-n^2_\pm)$, we can
apply \cref{lem:eigenvalues} \cite[App. D]{Siteur2014} to obtain information on the sign of the real parts of the respective eigenvalues:
\begin{lemma}\label{lem:eigenvalues}
    Let $\lambda_\pm =\alpha \pm \sqrt{\alpha^2 +\beta}$ for $\alpha, \beta\in\mathbb{R}\setminus\{0\}$. Then, for the real parts $\Re(\lambda_\pm)$,
    the following holds:
    \begin{equation*}
        \begin{array}{l|ll}
            & \beta > 0 & \beta < 0 \\
            \hline
            \alpha > 0
                & \Re(\lambda_+) > 0
                & \Re(\lambda_+) > 0 \\
            & \Re(\lambda_-) < 0
                & \Re(\lambda_-) > 0 \\
            \hline
            \alpha < 0
                & \Re(\lambda_+) > 0
                & \Re(\lambda_+) < 0 \\
            & \Re(\lambda_-) < 0
                & \Re(\lambda_-) < 0
        \end{array}
    \end{equation*}
\end{lemma}
We claim that regardless of $\alpha$, $(w_-, n_-)$ is unstable. We prove the claim by showing that $\beta > 0$ for $(w_-, n_-)$.
Indeed, we have that 
\begin{align*}
    w_- &= \frac{a+\sqrt{a^2-4m^2}}{2} > \frac{a}{2} > m
\end{align*}
and, recalling \cref{eq:w_eq}, we see that $n_- = \frac{m}{w_-} < 1 \implies \beta=m(1-n_-^2) > 0$, as $m>0$.

Next, we consider $(w_+,n_+)$.
We observe that 
\begin{equation}\label{eq:n>1}
    n_+ = \frac{a}{2m} + \frac{\sqrt{a^2-4m^2}}{2m} > 1+ \frac{\sqrt{a^2-4m^2}}{2m} > 1.
\end{equation}
Hence, $\beta = m(1-n_+^2) < 0$.
The stability of $(w_+,n_+)$ depends now on the sign of $\alpha$, as \cref{lem:eigenvalues} suggests that the equilibrium point is stable if and only if
\begin{equation}\label{eq:temp_klaus}
    \alpha < 0 \Leftrightarrow 1 + n_+^2 -m >0 \Leftrightarrow m < 1 + n_+^2.
\end{equation}
For $m<2$, we have that $1+n_+^2-m > 1 + n_+^2 -2 = n_+^2 -1 > 0$ by \cref{eq:n>1}. Hence, the equilibrium is asymptotically stable for $m<2 \land a>2m$.
The case $m\geq2$ is omitted here, since the magnitude of $m$ in real-world applications is usually below two \cite{Klausmeier1999,Siteur2014}.